\documentclass{amsart}

\usepackage{amsmath,amsthm,amssymb,amsfonts,eucal}
\usepackage{ifthen,epsfig,color,graphics}
\usepackage{rotating}

\newtheorem{thm}{Theorem}[section]
\newtheorem{prop}[thm]{Proposition}

\newtheorem{lem}[thm]{Lemma}

\theoremstyle{definition}
\newtheorem{defn}[thm]{Definition}

\newtheorem{alg}[thm]{Algorithm}

\newtheorem{rem}[thm]{Remark}
\newtheorem{step}{Step}

\theoremstyle{remark}

\renewcommand{\Re}{\text{Re}}
\renewcommand{\Im}{\text{Im}}
\renewcommand{\(}{\left(}
\renewcommand{\)}{\right)}  
\newcommand{\abs}[1]{\left\lvert #1\right\rvert}
\newcommand{\norm}[1]{\left\| #1\right\|}

\newcommand{\snorm}[1]{\norm{#1}}

\newcommand{\pnorm}[2]{\norm{#1}_{\rho, #2}}

\newcommand{\nbd}[2]{\mathcal{N}(#1, #2)}
\newcommand{\snbd}[2]{\mathcal{N}(#1, #2)}

\newcommand{\twovec}[2]{\genfrac{[}{]}{0pt}{}{#1}{#2}}
\newcommand{\field}[1]{\mathbb{#1}}
 
\newcommand{\FF}{\ensuremath{\field{F}}} 
\newcommand{\KK}{\ensuremath{\field{K}}} 
 
\newcommand{\IF}{\ensuremath{\field{I}\field{F}}} 
\newcommand{\IK}{\ensuremath{\field{I}\field{K}}} 
\newcommand{\CC}{\ensuremath{\field{C}}} 
\newcommand{\Ct}{\ensuremath{\field{C}^2}}

\newcommand{\RR}{\ensuremath{\field{R}}}
\newcommand{\Rt}{\ensuremath{\field{R}^2}}

\newcommand{\CPo}{\ensuremath{\field{C}\field{P}^1}}

\newcommand{\Henon}{H\'{e}non}
\newcommand{\Hypatia}{Hypatia}
\newcommand{\eps}{\epsilon}

\newcommand{\boxchcn}{box chain construction}
\newcommand{\boxcov}{\mathcal{B}}

\newcommand{\boxchmod}{box chain model} 
\newcommand{\boxchmods}{box chain models} 

\newcommand{\hyptwoalg}{Axis Metric Algorithm}

\newcommand{\showcomments}{no}

\newsavebox{\commentbox}
\newenvironment{comment}%
{\ifthenelse{\equal{\showcomments}{yes}}%
{\footnotemark
    \begin{lrbox}{\commentbox}
    \begin{minipage}[t]{1in}\raggedright\sffamily\small
    \footnotemark[\arabic{footnote}]}
{\begin{lrbox}{\commentbox}}}%
{\ifthenelse{\equal{\showcomments}{yes}}%
{\end{minipage}\end{lrbox}\marginpar{\usebox{\commentbox}}}
{\end{lrbox}}}

\newcommand{\drawfigtricolor}{\scalebox{.25}{\includegraphics{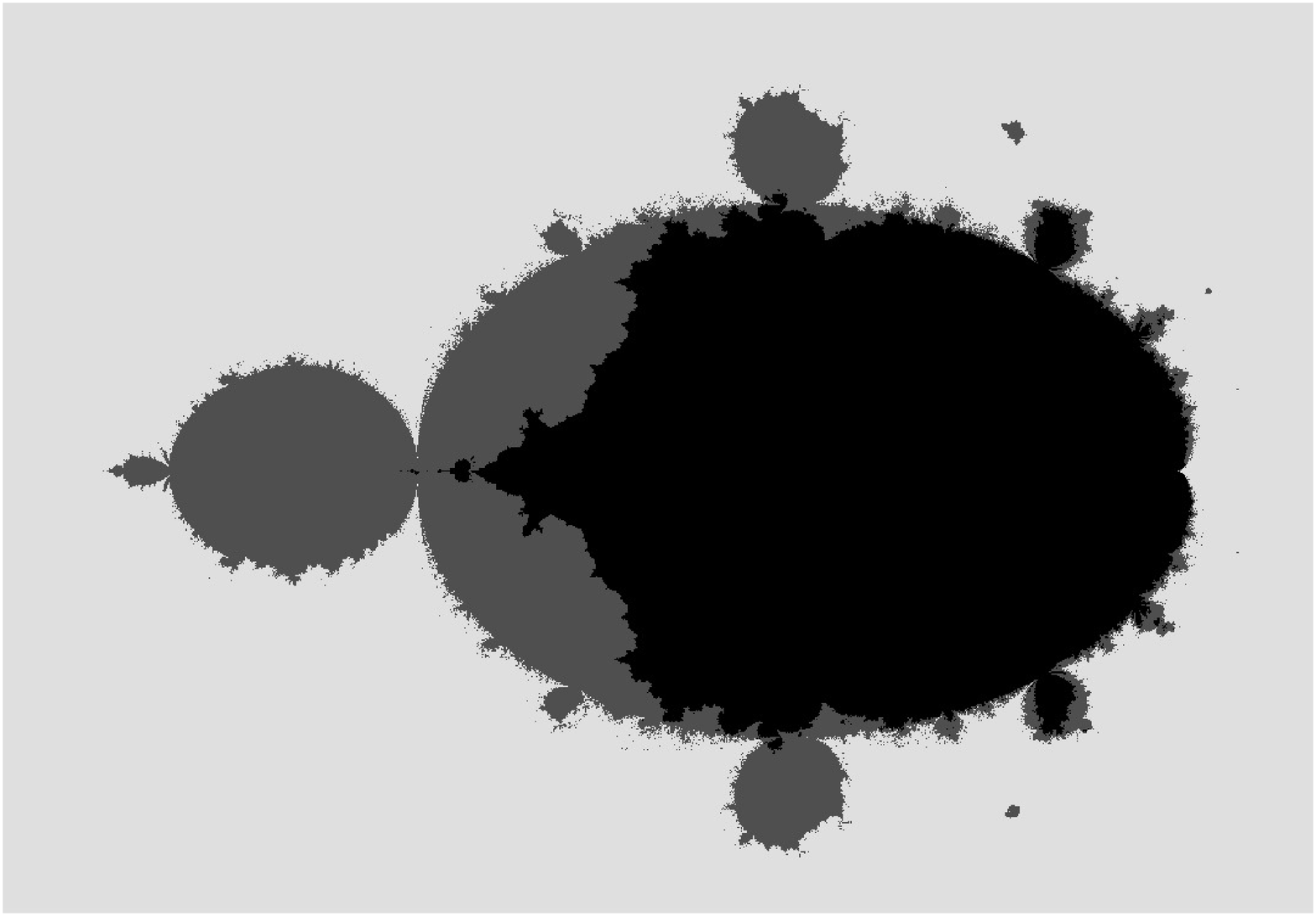}}}
\newcommand{\drawfigHenUMric}{\scalebox{.7}{\includegraphics{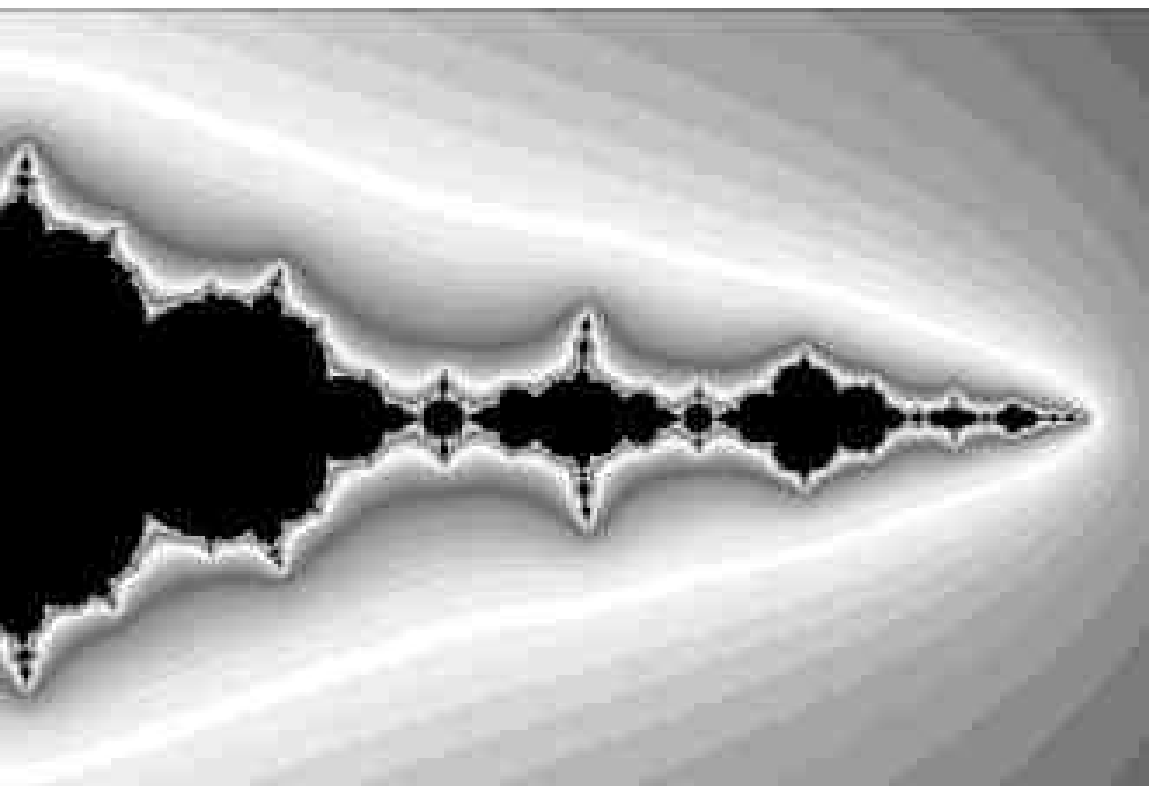}}}

\newcommand{\drawfighyptwoperoneA}{\scalebox{.4}{\includegraphics{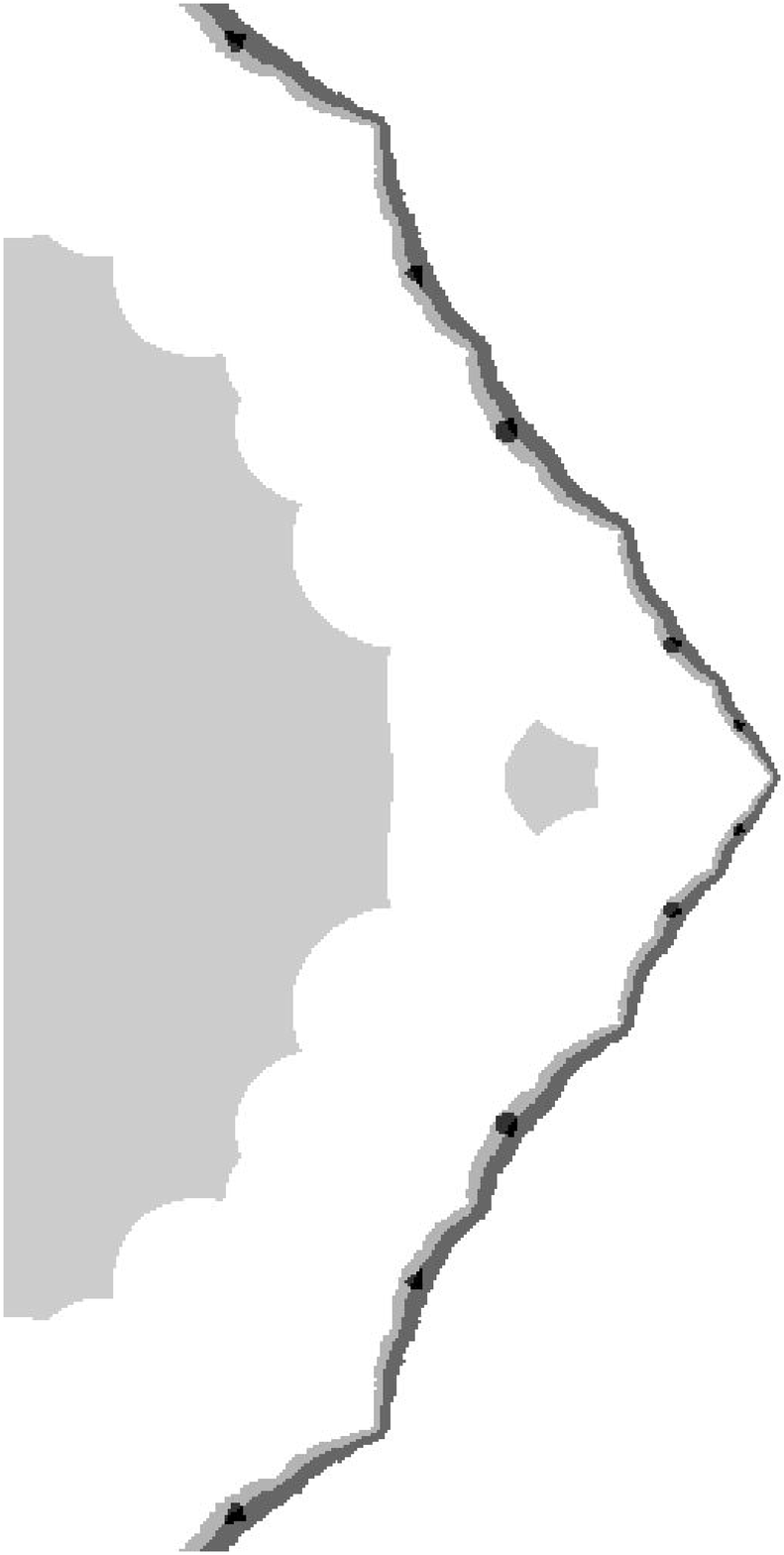}}}
\newcommand{\drawfighyptwocantorone}{\scalebox{.3}{\includegraphics{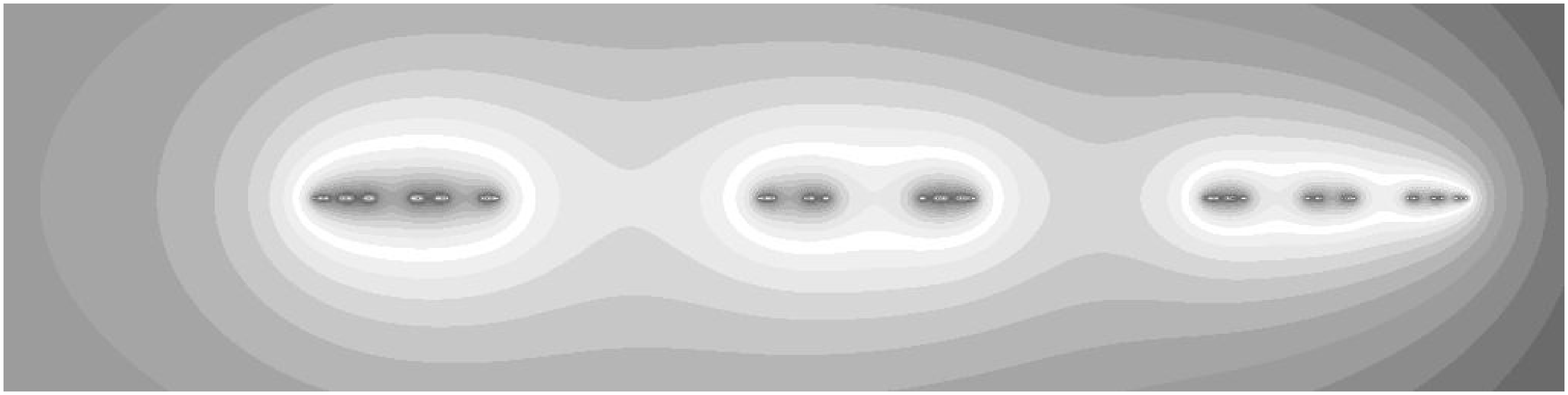}}}
\newcommand{\drawfighyptwocantortwo}{\scalebox{.5}{\includegraphics{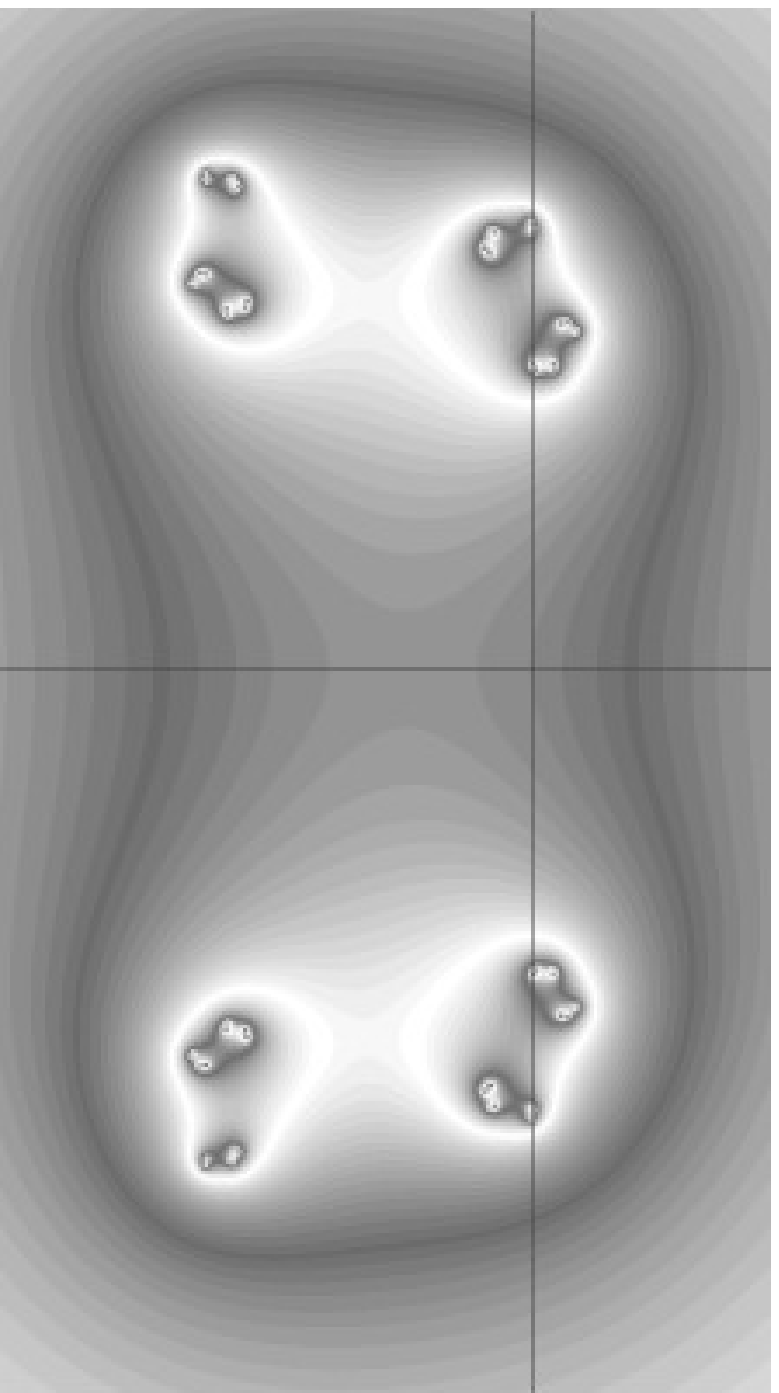}}}

\begin{document}

\title[Hyperbolicity of complex H\'{e}non maps]{A numerical method for constructing the hyperbolic structure 
of complex H\'{e}non mappings}
\author[S.L. ~Hruska]{Suzanne Lynch Hruska}
\thanks{Research supported in part by a grant from the National Science Foundation.}
\address{Department of Mathematics\\
Indiana University\\
Rawles Hall\\
Bloomington, IN 47405, USA
}
\email{shruska@msm.umr.edu}

\date{\today}
\begin{abstract}
For complex parameters $a,c$, we consider the H\'{e}non mapping $H_{a,c}: \Ct \rightarrow \Ct$,  given by $(x,y) \mapsto (x^2 +c -ay, x)$,  and its Julia set, $J$.
In this paper, we describe a rigorous computer program for attempting to construct a cone field in the tangent bundle
over $J$, which is preserved by $DH$, and a continuous norm in which $DH \ (\text{and } DH^{-1})$ uniformly expands the cones (and their complements).
We show a consequence of a successful construction is a proof
that $H$ is {hyperbolic} on $J$.  We give several new examples of hyperbolic maps,
 produced with our computer program, Hypatia, which implements our methods.
\end{abstract}



\renewcommand{\subjclassname}{\textup{2000} Mathematics Subject Classification}

\maketitle


\section{Introduction}
\label{sec:intro}


\begin{comment} All marginal comments will be removed for the submitted version.
\end{comment}
\begin{comment}
\textit{Point 1: henon maps are interesting to study (needs little justification)}.
\end{comment}
The {\Henon } family,  $H_{a,c} (x,y) = (x^2 +c -ay, x)$, 
has been extensively studied as a diffeomorphism of $\RR^2$,  with $a,c$ real parameters.
For example, Benedicks and Carleson show the existence of chaotic behavior
in the form of a strange attractor for some real {\Henon } maps in \cite{BC}.  Here we consider $H_{a,c}$ as a diffeomorphism of $\Ct$, and allow $a,c$ to be complex.
Foundational work on the dynamics of the complex {\Henon } family has been done by Bedford and Smillie  (\cite{BLS1, BLS2, BS6, BS9}), Hubbard (\cite{HOV1, HOV2, HPV}), and Fornaess and Sibony (\cite{FS1992c}).
However, basic questions
remain unanswered.

\begin{comment}
\textit{Point 2: there are methods available to study and understand hyperbolic henon mappings.}
\end{comment}
A natural class of maps to study are the hyperbolic maps, since hyperbolic maps 
generally have nontrivial (chaotic) dynamics, but are amenable to analysis.  
 A {\Henon } mapping $H_{a,c}$ is hyperbolic if its Julia set, $J_{a,c}$, is  a hyperbolic set for $H_{a,c}$.
(\textit{Hyperbolicity} and \textit{Julia set} are defined in Section 2.)
For {\Henon } mappings, hyperbolicity implies Axiom A,
which implies shadowing on $J$,
\textit{i.e.}, $\eps$-pseudo orbits are $\delta$-close to true orbits, and
structural stablity on $J$, \textit{i.e.}, in a neighborhood in
parameter space the dynamical behavior is of constant topological
conjugacy type.  
%
 %
Thus for a hyperbolic mapping, the dynamics on $J$ should be able to be
understood using combinatorial models.  
These properties make hyperbolic diffeomorphisms amenable
to exploration via
computers.

\begin{figure}
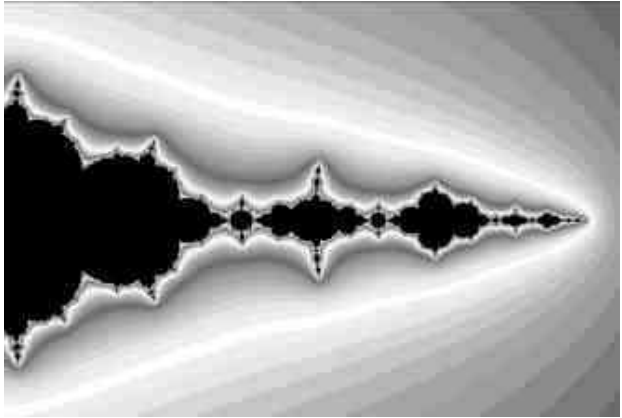

\begin{center}
  \drawfigHenUMric
\end{center}
\caption{\label{fig:HenUMric}  The filled Julia set for a {\Henon } mapping $H_{a,c}(x,y)$ $=$ $(x^2+c-ay,x)$, restricted to the 
unstable manifold of a saddle fixed point, with its natural 
parameterization.  
Here $a=.3, c=-1.17$, and the map has attracting cycles of  
periods one and three, which is impossible for quadratic polynomial maps of~$\CC$.
}
\end{figure}

\begin{figure}
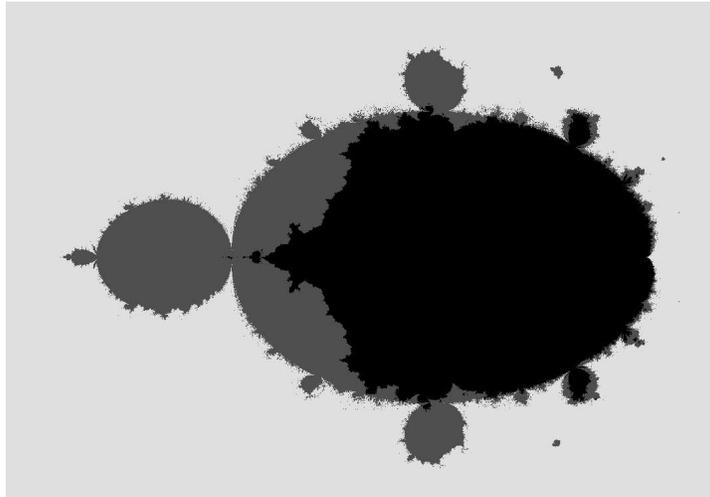

\begin{center}
\drawfigtricolor 
\end{center}
\caption{\label{fig:HenParam} A slice of {\Henon } parameter space, the 
$c$-plane, with $a=.3$. The innermost dark set
is an attempt to sketch the set of maps with
connected Julia sets, 
while the exterior is an attempt to sketch the maps with filled Julia set having empty interior 
(including complex horseshoes).
These regions have features reminiscent of the
parameter regions of one-dimensional polynomial maps.  This is an
intriguing parallel, suggesting that the {\Henon } parameter space may be
an equally rich arena of study.
}
\end{figure}

\begin{comment}
\textit{Point 3: Oliva and HP found interesting phenomena, and if the maps they studied were known to be hyp, one could try to prove these apparent phenomena are true.}
\end{comment}
%
Motivated by careful computer investigations,    
Oliva (\cite{OT}) provides a combinatorial model of the dynamics of some {\Henon } mappings, including for example, the mapping of Figure~\ref{fig:HenUMric}. The proposed model presupposes that the mapping is hyperbolic.  %
Hubbard and Papadantonakis (\cite{CUweb,HubKarl}) have more recently generated
pictures of slices of the {\Henon } parameter space, which attempt to sketch either the
 locus of maps with $J$ connected, or the locus of maps with $J$ having no interior (see  Figure~\ref{fig:HenParam}).
These and other computer investigations  
suggest that the dynamical behavior of the complex {\Henon } family is rich and subtle.  If certain of these mappings could be shown to be hyperbolic, then this would serve as a first step toward showing mathematically that these apparent phenomena actually occur.
%

\begin{comment}
\textit{Point 4: There are very few {\Henon } mappings known to be hyperbolic.}
\end{comment}
However, there are very few complex {\Henon } mappings known to be hyperbolic.  
Let us summarize what is known.
%
First, if
\begin{equation} \label{eqn:Horse}
\abs{c} > 2(1+\abs{a})^2,
\end{equation}
then $H_{a,c}|_J$ is
conjugate to the full 2-shift (so $J_{a,c}$ is a Cantor set), and the map is hyperbolic (compare Devaney and Nitecki
\cite{DN}, Oberste-Vorth \cite{OVT}, Morosawa, et al.\ \cite{MNTU}). 
In this case the mapping is called a {\em (complex) horseshoe}.  
The exterior in Figure~\ref{fig:HenParam} contains the set of horseshoes (among other types of maps).
Second, Hubbard and Oberste-Vorth (\cite{HOV2}) show that if
$P_c(x)=x^2+c$ is a hyperbolic polynomial, then there exists an $A(c)$ such
that if 
\begin{equation} \label{eqn:smallA}
0<~\abs{a}~< A(c),
\end{equation}
then $H_{a,c}|_{J}$ is topologically conjugate to the function induced by $P_c$
on the inverse limit ${\lim_{\leftarrow} (J, P_c)}$, hence $H_{a,c}$ is hyperbolic.  
%
Ishii and Smillie (\cite{IS}) have worked to obtain explicit estimates for the constant $A(c)$ in~(\ref{eqn:smallA}), but these estimates are relatively small.
%
%


\begin{comment}
\textit{Point 5: The purpose of this paper is to give an algorithm which can prove that {\Henon } mappings are hyperbolic}
\end{comment}
%
 Our broad goal is to develop computer algorithms 
with which we can rigorously describe the dynamics of any hyperbolic complex 
{\Henon } mapping.  In this paper we make a key step in that process, by developing a computer program which can establish if a complex {\Henon } mapping is hyperbolic, and if so, the program produces explicit information about how the map is hyperbolic; in particular, it builds two complementary cone fields in the tangent bundle over $J$ (the \textit{unstable} and \textit{stable} cones), and constructs a norm in which $DH \ (DH^{-1})$ preserves and unifomly expands the unstable (stable) cones.
%
%
%
%

Since hyperbolicity is structurally stable, a computer program with infinite resources should be able to prove hyperbolicity for any hyperbolic map.  However, non-hyperbolicity is an unstable condition, thus a computer program cannot be expected to recognize when a map is not hyperbolic.  (For example, $P(z) = z^{2} -1.5$  is presumably, but not provably, non-hyperbolic.)


This paper builds on the results of \cite{SLHone} (cf \cite{KMisch, Dell1, Osi, OsiCamp, Eiden}). There we describe an algorithm, called the \textit{\boxchcn}, which given $\eps >0$ finds a compact neighborhood, $\boxcov$, containing a $\delta= \delta(\eps)$-neighborhood of $J$,
and creates a finite graph, $\Gamma$, which models the $\eps$-dynamics of $H$ on $\boxcov$.
%
%
In this paper, we need only know the following about $\Gamma$.

\begin{defn} \label{defn:boxmodel}
Let $\Gamma = \Gamma (\mathcal{V}, \mathcal{E})$ be a directed graph,
with vertex set $\mathcal{V} = \{ B_k \}_{k=1}^N$  a finite collection of closed boxes in $\Ct$, having disjoint interiors, and  such that the union of the boxes $\mathcal{B} = \cup_{k=1}^N B_k$ contains $J$.  Suppose there is a $\delta>0$ such that 
$\Gamma$ contains an edge from $B_k$ to $B_j$ if the image $H(B_k)$ intersects a $\delta$-neighborhood of $B_j$, \textit{i.e.},  
\[
\mathcal{E} \supset \{ (k, j) \colon H(B_k) \cap \mathcal{N} (B_j,{\delta} ) \neq \emptyset \}.
\]
Further, assume $\Gamma$ is strongly connected, \textit{i.e.}, for each pair of vertices $B_k, B_j$, there is a path in $\Gamma$ from $B_k$ to $B_j$, and vice-versa.   Then we call $\Gamma$ a \textit{{\boxchmod } of $H$ on $J$}.
\end{defn}

In \cite{SLHone}, we describe how the {\boxchcn } builds strongly connected graphs modeling every basic set of $H$, for example, $J$ and any attracting periodic orbits.  In this paper, we let $\Gamma$ be the strongly connected graph component containing $J$, and we mostly ignore the others.

Our first task in this paper is to develop a discrete condition on $\Gamma$ we call \textit{box hyperbolicity},
 and show this condition implies hyperbolicity of $H$:
%
\begin{comment} Point 6 = Result (1): the discrete condition ``box hyp'' implies the continuous property  ``hyp.'' \end{comment}
\begin{thm} \label{thm:truehyptwo}
Let $\Gamma$ be a {\boxchmod } of $H$ on $J$.  If $\Gamma$ is box hyperbolic, then $H$ is hyperbolic on $J$. 
\end{thm}
Our definition of box hyperbolicity is inspired by our work in numerically establishing hyperbolicity in one complex dimension (see \cite{SLHtwo}). The difference is that in one dimension, hyperbolicity simply means expansion on $J$, whereas for {\Henon } mappings, hyperbolicity means a saddle property, \textit{i.e.}, expansion, contraction, and transversality.
The notion of box hyperbolicity is made precise in
 Section~\ref{sec:boxhyp}.
 Let us briefly describe this property here.  We begin with the cone field criterion for hyperbolicity.  
In a more general setting, Newhouse and Palis (\cite{New, NP1}) show that 
 an $f$-invariant set $\Lambda$  
is hyperbolic for $f$ iff there is a field of cones in the tangent bundle over $\Lambda$
such that
 $Df$ maps the cone field inside itself, and such that  in some norm,
$D f$ uniformly expands the cones, and $Df^{-1}$ uniformly expands the complements of the cones. Moreover, the field of cones need not be continuous in $x \in \Lambda$;
hence the cone field criterion for hyperbolicity yields a natural way to study the hyperbolic structure of a diffeomorphism using a computer (\textit{i.e.}, discretely). 
Here, we build cones which are constant on each box vertex of~$\Gamma$.   

To use the cone field criterion, we must find both a field of cones  preserved by $DH$ and a norm in which $DH$ expands the cones. 
 We cannot expect that $DH$ expands vectors in each cone with respect to the euclidean norm.  
 For example, in the euclidean norm,  
 $DH$ may be small for some cones over a pseudo-cycle, but larger in others, so that only the total cycle multiplier is larger than one.  Thus given a model $\Gamma$ of a map $H$, we attempt to build a discretized norm on the tangent bundle over $\Gamma$, which is designed to factor out the differences in $DH$ along cycles, so that in this new norm, $DH$ is expanding on every cone. Then we show that hyperbolicity in the discrete norm implies hyperbolicity for some continuous norm.

  

In order to test for box hyperbolicity on a given $\Gamma$, we have developed a computer algorithm we call the {\hyptwoalg}, designed to either prove that a given $\Gamma$ is box hyperbolic,
or describe which parts of $\Gamma$ are obstructions to proving box hyperbolicity.  
%
%
%
%
%
%
%
If a $\Gamma$ fails to be box hyperbolic, then either the map is not hyperbolic, or the boxes of $\Gamma$ are too large. Thus our approach is to attempt to prove box hyperbolicity on a sequence of models $\Gamma(n)$ with decreasing box size.  If $\Gamma(n-1)$ fails to be box hyperbolic, then to create $\Gamma(n)$ we could decrease the size of all the boxes, or use output of the {\hyptwoalg } on $\Gamma(n-1)$ in choosing which boxes need to be decreased in order to create a $\Gamma(n)$ more likely to be box hyperbolic.  If some $\Gamma(n)$ is found to be box hyperbolic, then $f$ is hyperbolic, and the program terminates. 
Thus a ``successful'' run of an implementation of our procedure gives a mathematical proof of hyperbolicity: if we can construct a $\Gamma$ which the {\hyptwoalg } shows to be box hyperbolic,
 then that mapping is
mathematically proven to be hyperbolic,  by Theorem~\ref{thm:truehyptwo}.

In designing the {\hyptwoalg } for verifying box hyperbolicity of some $\Gamma$, we build on the one-dimensional procedure described in \cite{SLHtwo}.  There we prove hyperbolicity of polynomial maps of $\CC$ by creating
a piecewise continuous (box constant) metric, under which the map is expanding on a neighborhood of  $J$.  To move up to two dimensions, and saddle-type hyperbolicity, we break down the problem 
into one dimensional pieces, then reassemble.
In particular, we build approximately invariant, box constant unstable and stable line fields, which will serve as axes for our cones.  Then we use the one dimensional algorithm twice, to
attempt to build a metric which is contracted on the stable directions, and
another metric which is expanded on the unstable directions.  These metrics
and line fields then determine a cone field with an induced metric, which we
use to test for box hyperbolicity.  
The {\hyptwoalg } is described in detail in Section~\ref{sec:hyptwoalg}.

\begin{comment}
\textit{Point 7: Result (2): For these 4 params, and run an implementation of hypatia, and found the maps to be hyp.}
\end{comment}
Finally, we have
 implemented our methods
 into a computer program called \textit{Hypatia}, and used the program to prove hyperbolicity of several {\Henon } mappings which were not previously known to be hyperbolic:





%
%

\begin{thm} \label{thm:examples}
The complex {\Henon } mappings, $H_{a,c}(x,y) =$ $(x^2+c-ay,x)$, with:
\[
(c,a) = (-0.3,0.1), (0, -0.22), (-3, -0.25), 
 \text{ and } (1.5, 0.5),
\]
are hyperbolic.
\end{thm}


Computer pictures suggest that the first two mappings of Theorem~\ref{thm:examples}, $H_{a,c}$ with $(c,a) = (-.3, .1)$ and $(c,a)=(0, -.22)$, are in the main cardioid,
with the only recurrent dynamics consist of a connected $J$ and one attracting fixed point, and that the latter two mappings are horseshoes, with
$J_{a,c}$ with $c=-3, a=-.25$ 
appearing to lie in $\Rt$, and $J_{a,c}$, with $c=1.5, a=.5$ appearing not to be contained in $\Rt$. 
Whether or not that is the case, each of the maps of Theorem~\ref{thm:examples} lies outside of the known regions defined by~(\ref{eqn:Horse}) and~(\ref{eqn:smallA}) with the Ishii-Smillie estimates.






All of the computations involved in proving Theorem~\ref{thm:examples}
were run on a Sun Enterprise E3500 server
with 4 processors, each $400$MHz UltraSPARC (though the multiprocessor was
not used) and $4$ GB of RAM.
\begin{footnote}{This server was purchased through an NSF SCREMS grant obtained
by the Department of Mathematics at Cornell University.}\end{footnote}
  When computations became overwhelming,
memory usage was the limiting factor.  
The C++, unix program, Hypatia, may be obtained from the author.

To conclude the introduction, we sketch the organization of the paper.
We give background on the dynamics of the {\Henon } family in Section~\ref{sec:bckgnd}.
In Section~\ref{sec:IA}, we briefly discuss  {\em interval arithmetic with directed rounding},
the method used to maintain rigor in our computer computations. 
In Section~\ref{sec:boxhyp}, we define box hyperbolicity, 
and we prove Theorem~\ref{thm:truehyptwo}, establishing that 
box hyperbolicity implies hyperbolicity.
In Section~\ref{sec:hyptwoalg}, we describe our computer procedure
for verifying box hyperbolicity, including the {\hyptwoalg}.  
Finally, in Section~\ref{sec:examples} we provide some data on
how we used  {\Hypatia } to prove hyperbolicity of each of the maps of Theorem~\ref{thm:examples}.

\bigskip
\textbf{Acknowledgements.}
These results were primarily accomplished as my PhD thesis at Cornell University (\cite{SLHT}). I am grateful to John Smillie for providing guidance on the project, John Hubbard for inspiration, Greg Buzzard, and Warwick Tucker  for many helpful conversations,
Eric Bedford, James Yorke, and John Milnor for advice on the 
preparation of this paper, and Robert Terrell for technical support. I would also like to thank the referees and editors for providing comments which I helped me to significantly increase the clarity of the presentation.

\section{Background}
\label{sec:bckgnd}

\subsection{The {\Henon } Family}
\label{sec:henonintro}

Polynomial diffeomorphisms of $\Ct$ necessarily have polynomial inverses,
thus are often called polynomial automorphisms.
Friedland and Milnor (\cite{FM}) showed that
polynomial automorphisms of $\Ct$ break down into two categories.
\textit{Elementary} automorphisms have simple dynamics, and are
polynomially conjugate to a diffeomorphism of the form $(x,y) \mapsto
(ax+b, cy+p(x))$ ($p$ polynomial, $a,c \neq 0$). \textit{Nonelementary}
automorphisms are all conjugate to
finite compositions of \textit{generalized {\Henon }
mappings}, which are of the form $f(x,y) = (p(x)-ay,x)$, where $p(x)$ is a
monic polynomial of degree $d>1$ and $a \neq 0$.

To clarify the situation, one can define a \textit{dynamical degree} of a 
polynomial automorphism of $\Ct$.  If \textit{deg}$(f)$ is the maximum of the 
degrees of the coordinate functions, the dynamical degree is 
\[ 
d = d(f) = \lim_{n \to \infty} (\textit{deg}(f^{n}))^{1/n}. 
\] 
This degree is a conjugacy invariant.  Elementary automorphisms have 
dynamical degree $d=1$.  A nonelementary automorphism is conjugate to some
automorphism whose polynomial degree is equal to its dynamical degree. 
Without loss of generality, we assume such $f$ are finite
compositions of generalized {\Henon } mappings, rather
than merely conjugate to mappings of this form.

Thus, the quadratic, complex {\Henon } family $H_{a,c} (x,y)=(x^2+c-ay,x)$
represents the dynamical behavior of the simplest class of nonelementary
polynomial automorphisms; those of dynamical degree two.  In this paper, we usually
use the letter $f$ for a polynomial diffeomorphism of $\Ct$ with
$d(f)>1$, and $H$ for a (degree two) {\Henon } mapping.  We state results in Section~\ref{sec:boxhyp}
for the more general $f$, but in explaining the procedure for verifying hyperbolicity in Section~\ref{sec:hyptwoalg}, we concentrate on the case of $H_{a,c}$.

\subsection{Drawing Meaningful Pictures}
\label{sec:henwuppic}

For a polynomial map $P$ of $\CC$, the {\em filled Julia 
set}, $K$, is the set of points
whose orbits are bounded under $P$;  the {\em Julia set},
$J$,  is the topological boundary of $K$. 
For a polynomial diffeomorphism $f$, like $H_{a,c}$, there are 
corresponding Julia sets:
\begin{itemize}
\item $K^+ (K^-)$ is the set of points whose orbits are bounded under
$f (f^{-1})$ \\
and $K = K^+ \cap K^-$ is called the \textit{filled Julia set};
\item $J^{\pm} = \partial K^{\pm}$ (topological boundary) \\
and $J = J^+ \cap J^-$ is called the \textit{Julia set}.
\end{itemize}

Filled Julia sets are the (chaotic) invariant sets which can be easily sketched by
computer, on any two-dimensional slice.  Hubbard has suggested the
following method for drawing a dynamically significant slice of the Julia
set, by parameterizing an unstable manifold.  This method has been
implemented by Karl Papadantonakis in FractalAsm (\cite{HubKarl, CUweb}).  Figures~\ref{fig:HenUMric}, \ref{fig:cantorone}, and~\ref{fig:cantortwo} were generated using FracalAsm.

Let $f$ be a diffeomorphism of $\Ct$.  If $p$ is a periodic point of 
period $m$, and the eigenvalues
$\lambda, \mu$ of $D_pf^m$ satisfy $\abs{\lambda} > 1 > \abs{\mu}$ (or
vice-versa), then $p$ is a {\em saddle periodic point}.
The large (small) eigenvalue 
is called the unstable
(stable) eigenvalue.
If $p$ is a saddle periodic point,
then the {\em stable manifold of $p$} is
$ 
W^s(p) = \{q \colon d(f^n(q),f^n(p)) \rightarrow 0 \text{ as }
n\rightarrow \infty\}, 
$
and the {\em unstable manifold of $p$} is
$ 
W^u(p) = \{q \colon d(f^{-n}(q),f^{-n}(p)) \rightarrow 0 \text{ as }
n\rightarrow \infty\}. 
$
If $p$ a saddle periodic point of $f$, then $W^u(p) \ (W^s(p))$ is biholomorphically 
equivalent 
to $\CC$, and on $W^u(p) \ (W^s(p))$, $f$ is conjugate 
to multiplication by the unstable
(stable) eigenvalue of $D_pf$.

When $\abs{a} \neq 1$, except on the curve of equation $4c = (1+a)^2$,
the map $H_{a,c}$ has at least one saddle fixed point, $p$, 
(\cite{HubKarl}). The unstable manifold 
$W^u(p)$ has a natural parametrization $\gamma \colon \CC \to W^u(p)$ 
given by
\[
\gamma(z) = 
\lim_{m\to\infty} \gamma_m(z) = 
\lim_{m\to\infty} 
H^m \( p + \frac{z}{\lambda_1^m}\mathbf{v_1} \),
\]
where $\lambda_1$ is the unstable eigenvalue of $D_pH$ and
$\mathbf{v_1}$ is the associated eigenvector.
This parametrization has the property that
$
H(\gamma(z)) = \gamma(\lambda z),
$
and any two parametrizations with this property differ by scaling the 
argument.

Observe that since $W^u(p) \subset K^-$, to get a picture of $K$ in
$W^u(p)$ we need only color pixels black which are guessed to be in $K^+$.
To sketch the picture, we approximate $\gamma$ by say $g=\gamma_{100}$ in
a region in the plane, $B=\{z=x+iy \colon a \leq x \leq b, c \leq y \leq
d\}$.  Then an escape threshold is chosen, like $10$, and then 
for each $z\in B$, if $\norm{H^{n}(g(z))}  < 10$ for all $n<100$, we
say $g(z)\in K^+$ and color it black.  Otherwise, color according to which
iterate $\norm{H^n(g(z))}$ first surpassed $10$.

\subsection{Hyperbolicity}
\label{sec:hypintro}

In this subsection, let $f$ be a diffeomorphism of a manifold $M$, and let $\Lambda$ be a compact, $f$-invariant set.  First we recall the standard definition of hyperbolicity (see \cite{Rob}).

\begin{defn} \label{defn:hyp}
$\Lambda$
is {\em hyperbolic} for $f$ if at each $x$ in $\Lambda$, there is a splitting of the tangent bundle $T_xM
= E^s_x \oplus E^u_x$, which varies
continuously with $x \in \Lambda$, such that:
\begin{enumerate}
 \item $f$ preserves the splitting, \textit{i.e.}, $D_xf(E^s_x) = 
E^s_{fx}$, and $D_xf(E^u_x) = E^u_{fx}$, and 
\item $Df \ (Df^{-1})$ expands on $E^u \ (E^s)$ uniformly, \textit{i.e.}, there exists a constant $\lambda > 1$ and a norm $\norm{\cdot}_x$ on $T_{\Lambda} M$, continuous for $x \in \Lambda$, for which
\begin{eqnarray*}
 \norm{D_xf(\mathbf{w^u})}_{f(x)} & \geq & \lambda \norm{\mathbf{w^u}}_x, \text{ for } \mathbf{w^u} \in E^u_x, \text{ and }\\
 \norm{D_xf^{-1}(\mathbf{w^s})}_{f^{-1}(x)} & \geq & \lambda \norm{\mathbf{w^s}}_x, \text{ for } \mathbf{w^s} \in E^s_x.
\end{eqnarray*}
\end{enumerate}
\end{defn}

As noted in the Introduction, Newhouse and Palis (\cite{New, NP1}) show hyperbolicity can be described using a {\em cone field}.  To define a cone $C_x$ for each point $x$ in $\Lambda$, we need a
splitting $T_x M = E_{1x} \oplus E_{2x}$, and a
positive real-valued function $\eps(x)$ on~$M$.  Then define the
$\eps(x)$-sector $S_{\eps(x)}(E_{1x}, E_{2x})$ by
\[ 
S_{\eps(x)}(E_{1x}, E_{2x}) 
= \{ (\mathbf{v}_1, \mathbf{v}_2) \in E_{1x} \oplus E_{2x}:
      \norm{\mathbf{v}_2} \leq \eps(x)\norm{\mathbf{v}_1}  \}.
\]
Then $C_x = S_{\eps(x)}$.
Newhouse and Palis show that $\Lambda$  
is hyperbolic for $f$ iff there is a field of cones $\{ C_x \subset T_xM \colon x \in \Lambda\}$, 
a constant $\lambda > 1$, and a continuous norm 
$\norm{\cdot}$, such that
 $Df$ preserves the cone field, \textit{i.e.}, $D_x f(C_x) \subset C_{f x}$, and such that in this norm,
$D_x f \ (D_x f^{-1})$ uniformly expands vectors in $C_x \ (T_x M - C_x)$;
moreover, the field of cones $x \rightarrow C_x$ need not be continuous.
Our computer algorithm for verifying hyperbolicity actually combines these two notions, as we will see in Section~\ref{sec:boxhyp}.



Bedford and Smillie (\cite{BS2}) have shown that for $f$ a polynomial 
diffeomorphism of $\Ct$, with $d(f)>1$,  $f$ is hyperbolic  on its Julia set, 
$J$,
iff $f$ is hyperbolic on its chain recurrent set, $\mathcal R$, iff
$f$ is hyperbolic on its nonwandering set, $\Omega$.
Thus we say $f$ is hyperbolic if any of these conditions holds.
In fact, in \cite{BS1} Bedford and Smillie show that if $f$ is hyperbolic, then
$\mathcal R$ and  $\Omega$ are both equal to $J$ union
 finitely many attracting periodic  orbits.
Thus for hyperbolic polynomial diffeomorphisms of $\Ct$, 
the 
basic sets are $J$ and the attracting periodic orbits.

\section{Interval arithmetic}
\label{sec:IA}
\bigskip

In order to genuinely prove dynamical properties, we use in 
Hypatia a method of controlling round-off error in the computations,
called {\em interval arithmetic with directed rounding} (IA).  This method was recommended
by Warwick Tucker, who used it in his recent computer proof that the
Lorenz differential equation has the conjectured geometry (\cite{War}).

In fact we use IA not only to control error, but we take advantage of the
structure of this method in our algorithms and implementation.  
 We thus give a very brief description of IA below, and refer
the interested reader to \cite{MooreIA1, MooreIA2, GenIA}.

On a computer, we cannot work with real numbers, rather we work over
the finite space $\FF$ of numbers representable by binary floating point
numbers no longer than a certain length.  For example, since the number 
$0.1$ is not a dyadic rational, it has an infinite binary expansion.  
The computer cannot encode  this number exactly.  Instead,
the basic objects of arithmetic are not real
numbers, but closed intervals, $[a,b]$, with end points in some fixed field $\KK$.  
Denote this space of intervals by $\IK$.  
The operation of addition of two 
intervals $[a,b],$ $[c,d] \in \IK$ is defined by:
$
[a,b] + [c,d] = \left[  a+c,
b+d  \right].
$
Hence if $x\in [a,b]$ and $y\in [c,d]$, then $x+y \in [a,b]+[c,d]$.

The other operations are defined analogously, for example:
\begin{align*}
[ a, b ] - [ c, d ]  & = 
\left[  a-d , b-c  \right], \text{ and } \\
[ a, b ] \times [ c, d ] & = 
\left[  \min(ac, ad, bc, bd), \max(ac, ad, bc, bd) \right] 
\end{align*}

However, an arithmetical operation on two numbers in $\FF$ may not have a result in $\FF$. Thus to implement rigorous IA we use the idea of {\em directed rounding} to round outward the result of any operation.  For example,
\[
[a,b] + [c,d] = \left[ \left\downarrow a+c \right\downarrow, \left\uparrow 
b+d \right\uparrow \right],
\]
 where $\left\downarrow x \right\downarrow$ is the largest number in $\FF$  
strictly less than $x$ (\textit{i.e.}, $x$ rounded down), and 
$\left\uparrow x 
\right\uparrow$ 
is the smallest number in $\FF$ strictly greater than $x$ 
(\textit{i.e.}, $x$ 
rounded up).  For any $x\in \RR$, let Hull$(x)$ be the smallest interval in $\IF$ which contains $x$.  Thus, if $x \in \FF$, then Hull$(x) = [x,x]$.  If $x \in \RR \setminus \FF$, then Hull$(x) = [\left\downarrow x \right\downarrow, \left\uparrow x \right\uparrow]$.

Thus, if the user is interested in a computation involving real numbers,
then IA with directed rounding performs the computation using intervals in $\IF$ which 
contain those real numbers, and gives the answer as an interval 
in $\IF$ which contains the real answer.
In higher dimensions, IA operations can be carried out component-wise, on 
{\em interval vectors}.  

One must think carefully
about how to use IA in each arithmetical calculation. 
For example, it can create problems by
propagating increasingly large error bounds.
Iterating a polynomial diffeomorphism like  
$H_{a,c}$ on an interval vector which is not very close to an 
attracting period cycle 
will give a tremendously large interval vector after only a few iterates.  
That is, suppose 
$B=[a,b]\times [c,d]$ is an interval vector in $\CC$, and one attempts to 
compute a box containing $P_c^{10}(B)$, for $P_c(z) = z^2+c$, by:
\begin{tabbing}
\hspace*{.25in} \= \hspace{.25in} \= \hspace{.25in} \kill
\>for $j$ from $1$ to $10$ do \\
\>\>$B = P_c(B)$
\end{tabbing}
then the box $B$ will likely grow so large that its defining bounds become 
machine $\infty$, \textit{i.e.}, the largest floating point in $\FF$.
Similarly, one would also never want to try to compute
$
D_{B_n}H \circ \cdots \circ D_{B_1}H \circ D_{B_0}H (\mathbf{u}),
$
for a vector $\mathbf{u} \in \Ct$, since the entries would blow up (see Algorithm~\ref{alg:hyptwoalg}, Step~\ref{alg:USdirs}).

Our construction involving boxes in $\Ct$ as the basic numerical objects is designed to be efficiently manipulated with IA. 
For all of our rigorous computations, we use
IA routines provided by the PROFIL/BIAS package, available at 
\cite{PBIA}.

\section{Characterizing box hyperbolicity}
\label{sec:boxhyp}

In this section we 
define box hyperbolicity for a {\boxchmod } $\Gamma$, 
in Definition~\ref{defn:boxhyp}, and show that if
$\Gamma$ is box hyperbolic,
then $f$ 
satisfies the standard definition of hyperbolicity, proving Theorem~\ref{thm:truehyptwo}.
Further, we show in 
Proposition~\ref{prop:conepresexp} that box hyperbolicity is equivalent to a simple condition in linear algebra.  Throughout this section, assume $f$ is a polynomial diffeomorphism of $\Ct$, with $d(f) > 1$, and let $\Gamma = (\mathcal{V}, \mathcal{E})$ be a {\boxchmod } of $f$ on $J$.

\begin{defn} \label{defn:boxhyp}
%
Suppose for each box $B_k$ in $\mathcal{V}$,  we have some nondegenerate indefinite Hermitian form, $Q_k \colon \Ct \to \RR$. 
Define ${C}^u_k := \{ \mathbf{w} \colon Q_k(\mathbf{w}) \geq 0 \}$, as the {\em unstable cones}, and define their complements as the {\em stable cones}:  ${C}^s_k := \Ct \setminus {C}^u_k = \{ \mathbf{v} \colon Q_k(\mathbf{v}) < 0 \}$.   

We say that $\Gamma$ is {\em box hyperbolic} if 
$Df$ ($Df^{-1}$) preserves and expands the unstable (stable) cones, with respect to $\{Q_k\}$, \textit{i.e.}, 
for every edge $(k,j) \in 
\mathcal{E}(\Gamma)$, and every $z \in B_k$:
\begin{enumerate}
\item if $\mathbf{w} \in  C^u_k$, then 
    $D_{z}f(\mathbf{w}) \in C^u_j$ and
    ${Q_j(D_{z}f(\mathbf{w}))} \ > \ 
    {Q_k(\mathbf{w})};$
    
\item if $\mathbf{v} \in  C^s_j$, then
    $\left[D_{z}f\right]^{-1} (\mathbf{v}) \in  C^s_k$ and
     $-{Q_k\(\left[D_{z}f\right]^{-1}(\mathbf{v})\)} \ > \ 
     -{Q_j(\mathbf{v})}$.
\end{enumerate}     

\end{defn}

In fact, a given cone $C^u_k$ determines an associated Hermitian form $Q_k$ up to scaling.  Finding an appropriate choice of scaling for each $Q_k$ is how we determine a metric for which $H$ is hyperbolic.  To prove Theorem~\ref{thm:truehyptwo}, we first use a partition of unity argument to smooth out the forms $\{ Q_k \}_{k=1}^N$ into a continuous field of forms $\{ Q_z \}_{z \in J}$ (Definition~\ref{defn:Qz}), then define a riemannian metric induced by  $\{ Q_z \}_{z \in J}$ (Definition~\ref{defn:pnorm}), and show that $H$ is hyperbolic on $J$ in this new riemannian metric.


\subsection{Box hyperbolicity implies hyperbolicity}

Here our goal is to prove  Theorem~\ref{thm:truehyptwo}, that if $\Gamma$ is box hyperbolic, then $f$ is hyperbolic on $J$, as in Definition~\ref{defn:hyp}.  Part of the proof is
very similar to the one dimensional analog, proved in \cite{SLHtwo}, 
in that we use a partition of unity
to smooth out a discrete norm. But before we deal with the norm, we verify that box hyperbolicity implies the existence of a continuous splitting preserved by the map.

\begin{lem} \label{lem:splitting}
If $\Gamma$ is box hyperbolic,
 then there exists
a splitting of the tangent bundle $T_{z}\Ct = E^s_z
\oplus E^u_z$, for each $z$ in $J$, which varies
continuously with $z$ in $J$,  such that $f$ preserves the splitting, 
\textit{i.e.}, $D_zf(E^s_z) = E^s_{fz}$, and
$D_zf(E^u_z) = E^u_{fz}$.  Further, for each $z\in B_k$, 
$E^u_z \subsetneq {C}^u_k$, and 
$E^s_z \subsetneq {C}^s_k$.
\end{lem}

\begin{proof}
Recall that Newhouse and Palis show that a
diffeomorphism $f$ is hyperbolic if there is a field of cones
${C}_z$ (not necessarily continuous) which is preserved and
expanded by $Df$, such that the complements are expanded by $Df^{-1}$.  
In their proof (\cite{NP1}), they first show that the existence of a
cone field preserved by $Df$ implies the existence of a continuous
splitting preserved by $f$, with the unstable (stable) directions lying
inside the unstable (stable)  cones. Box hyperbolicity gives a cone field
preserved by $Df$. Thus we have cones ${C}_z = {C}^u_k$,
if $z$ is in box $B_k$ (make some consistent choice of box $B_k$
containing $z$, for the benefit of points on the boundaries of the closed
boxes).  Thus by the proof in \cite{NP1}, we have the existence of the
continuous splitting preserved by~$Df$.
\end{proof}


It is more natural for computer calculations to use the $L^{\infty}$ metric on $\RR^{2n} = \CC^n$, rather than euclidean.  Thus, throughout the rest of this paper, $\norm{\cdot}$ will denote this norm, \textit{i.e.}, if $x = (x_1, \ldots, x_n)$, then
\[
\norm{x} = \max \{  \abs{\Re(x_k)}, \abs{\Im(x_k)} \colon 1 \leq k \leq n\}.
\]
If $B$ is a set of points in $\Ct$, we denote by $\snbd{B}{\delta}$ the 
$\delta$-neighborhood of the set $B$ in this metric.

When we say {\em box}, we mean a ball about a point in this norm.  Thus a box is also a vector of intervals, so boxes are neighborhoods which are easily manipulated with interval arithmetic.

We will need to measure the angle between pairs of lines in $\Ct$,
like $E^u_z$ and $E^s_z$, or $E^u_z$ and $E^u_x$.  To do so, we view the
set of lines through the origin in $\Ct$ as the projective space $\CPo =
\hat{\CC}$.  Then the spherical metric on $\CPo$ induces the following
metric.

\begin{defn} \label{defn:spheremet}
If $\mathbf{v} = \twovec{v_1}{v_2}, \mathbf{w} = \twovec{w_1}{w_2}$ are
vectors in $\Ct$, define the distance between the directions 
$\mathbf{\CC v}, \mathbf{\CC w}$ to be
\[
\sigma(\mathbf{v},\mathbf{w}) = \sin^{-1} \( 
	\frac{\abs{v_1 w_2 - v_2 w_1}}{\norm{\mathbf{v}} \norm{\mathbf{w}}} \).
\]
\end{defn}
Note that $\sigma(\mathbf{v},\mathbf{w}) \in [0,\pi/2]$, and for any complex
numbers $\alpha, \beta$,
\[
\sigma(\mathbf{v},\mathbf{w}) = \sigma(\mathbf{\alpha v},\mathbf{\beta w}).
\]
Think of this metric as measuring the angle 
between the complex lines.
If $\mathbf{v}$ is a vector and $\mathbf{W}$ is a collection of vectors in $\Ct$, we define $\sigma(\mathbf{v}, \mathbf{W}) = \inf\{\sigma(\mathbf{v},\mathbf{w}) \colon \mathbf{w} \in \mathbf{W} \}$. 

In the next two lemmas, we quantify how, in the metric $\sigma$, the unstable (stable) lines from the splitting of Lemma~\ref{lem:splitting} are strictly inside the unstable (stable) cones.

\begin{lem} \label{lem:EzinCk}
Let $\Gamma$ be box hyperbolic.
Then there exist $d_1 >0$ and $\delta_1 >0$ such that if $z\in J \cap 
\snbd{B_k}{\delta_1}$, then $\sigma(E^u_z, {C}^s_k) \geq d_1$ and
$\sigma(E^s_z, {C}^u_k) \geq d_1$.
\end{lem}

\begin{proof}
First, note that by compactness of $J$ and the fact that the line fields 
are contained in the interior of the cones, there exists a $d_0>0$ such 
that
\[
d_0 \leq \min \{ \sigma(E^u_z, {C}^s_k) \colon z\in J\cap B_k  \} 
\text{ \ and \ }
d_0 \leq \min \{ \sigma(E^s_z, {C}^u_k) \colon z\in J\cap B_k \}.
\]

Let $d_1 = d_0/2$.  
Next, by compactness of $J$ and continuity of the 
splitting, there exists a $\delta_1 >0$ such that for any $x,z \in J$
with $\snorm{x-z} < \delta_1$, we have $\sigma(E^u_z, E^u_x) < d_1$
and $\sigma(E^s_z, E^s_x) < d_1$.

Now let $z\in J \cap \snbd{B_k}{\delta_1}$.  Since $z$ is not necessarily 
in $B_k$, let $m$ be such that $z\in B_m$, and $x$ 
is a point satisfying $x\in J\cap B_m \cap B_k$ and $\snorm{x-z} < 
\delta_1$.  Then $\sigma(E^u_z, E^u_x) < d_1$ and $\sigma(E^s_z, E^s_x) < 
d_1$.
Since $x\in B_k$ we have $\sigma(E^u_x, 
{C}^s_k) \leq d_0$ and $\sigma(E^s_x, {C}^u_k) \leq d_0$.
Hence, $\sigma(E^u_z, {C}^s_k) \geq d_0 - d_1 = d_1$, and
$\sigma(E^s_z, {C}^u_k) \geq d_1$.
\end{proof}

\begin{lem} \label{lem:EfzinCj}
Let $\Gamma$ be box hyperbolic.
If $B_k, B_j \in \mathcal{V}$ and $z \in J$
satisfies $z\in \snbd{B_k}{\delta_1}$ and $f(z) \in \snbd{B_j}{\delta_1}$,
then $\sigma(E^u_{fz}, {C}^s_j) \geq d_1$ and 
$\sigma(E^s_{fz}, {C}^u_j) \geq d_1$.
\end{lem}

\begin{proof}
This lemma follows directly from Lemma~\ref{lem:EzinCk} applied to $f(z)$ 
instead of~$z$.
\end{proof}

Before the next step, we need a lemma from \cite{SLHone}.   

\begin{lem}[\cite{SLHone}] \label{lem:imageboxsize}
There exists an $\eta>0$ so that if $B_k, B_j \in \mathcal{V}$, with $z \in \nbd{B_k}{\eta}$ and $f(z) \in \nbd{B_j}{\eta}$, then there is an edge from $B_k$ to $B_j$ in $\Gamma$.
\end{lem}

To prove this lemma, we used the assumption that $f$ was a polynomial mapping of degree $d > 1$, and the fact that by Definition~\ref{defn:boxmodel}, there is a $\delta>0$ such that 
there is an edge from $B_k$ to $B_j$ 
if a $\delta$-neighborhood of $f(B_k)$ intersects $B_j$.  Now, we get:

\begin{lem}  \label{lem:deltaQexp}
Let $\Gamma$ be box hyperbolic.
Then there is a $\tau >0$ such that for any  $B_k, B_j \in 
\mathcal{V}$ and any $z \in J$ such that $z\in \snbd{B_k}{\tau }$ and 
$f(z) \in \snbd{B_j}{\tau }$, we have
\begin{enumerate}
\item if $\mathbf{w}\in E^u_z$, then 
	$Q_j\( D_zf(\mathbf{w}) \) > Q_k (\mathbf{w})$;
\item if $\mathbf{v}\in E^s_{fz}$, then
  	 $-{Q_k\(\left[D_zf\right]^{-1}(\mathbf{v})\)} \ > \
     -{Q_j(\mathbf{v})}.$
\end{enumerate}
\end{lem}

\begin{proof}
Among additional requirements given below, let $\tau $ be less 
than $\eta$ from 
Lemma~\ref{lem:imageboxsize}.
Then for any $z\in J$ such that $z \in  \snbd{B_k}{\tau }$ and $f(z) \in 
\snbd{B_j}{\tau }$, there is an edge in $\Gamma$ from $B_k$ to $B_j$,
\textit{i.e.}, $(k,j) \in \mathcal{E}$.

Since $J$ is compact, and by Lemmas~\ref{lem:EzinCk} 
and~\ref{lem:EfzinCj}, there exists $d_2\geq 0$ such that:
\[
d_2 \leq \min\{ Q_j(D_xf(\mathbf{u}_x)) - Q_k(\mathbf{u}_x) 
\colon x\in B_k, (k,j) \in \mathcal{E},
\mathbf{u}_x \in E^u_x, \norm{\mathbf{u}_x} =1 \}.
\]

Let $\eps = d_2/3$.  By continuity of $D_xf$ and the splitting, 
there is a $\tau < \eta$ so that
for any $x,z \in J$ with $\norm{x-z} < \tau $,  $z \in  
\snbd{B_k}{\tau }$, and $f(z) \in
\snbd{B_j}{\tau }$, we have
\begin{eqnarray*}
\abs{Q_k(\mathbf{u}_z)-Q_k(\mathbf{u}_x)} &<& \eps, \\
\abs{Q_j(D_xf(\mathbf{u}_x))-Q_j(D_xf(\mathbf{u}_z))} &<& \eps,
 \text{ and } 
\\
\abs{Q_j(D_xf(\mathbf{u}_z))-Q_j(D_zf(\mathbf{u}_z))} &<& \eps.
\end{eqnarray*}
Then $Q_j(D_zf(\mathbf{u}_z)) - Q_k(\mathbf{u}_z) \geq d_2 - 3\eps >0$.

Now since $Q(\mathbf{w})$ is a Hermitian form,
$Q(\alpha \mathbf{w}) = \abs{\alpha}^2 Q(\mathbf{w})$
for any $\alpha\in \CC$.  Thus by linearity of $Df$, the above
result for $\mathbf{u}_z$
implies the same result for any $\mathbf{w} \subset E^u_z$. Hence we have
Condition~$(1)$.

The proof of $(2)$ is analogous.  Let $d_3>0$ satisfy:
\[
d_3 \leq \min \{ Q_j(\mathbf{s}_{fz}) - Q_k(D_zf^{-1}(\mathbf{s}_{fz})) 
\colon x\in B_k, (k,j) \in \mathcal{E}, \mathbf{s}_{fz}\in E^s_{fz},
\norm{\mathbf{s}_{fz}} =1 \}.
\]
Let $\eps' = d_3/3$.  Then further restrict $\tau $ so that 
for any $x,z \in J$ with $\norm{x-z} <\tau $,  $z \in
\snbd{B_k}{\tau }$, and $f(z) \in
\snbd{B_j}{\tau }$, we have
\begin{eqnarray*}
\abs{Q_j(\mathbf{s}_{fz}) - Q_j(\mathbf{s}_{fx})} &<& \eps', \\
\abs{Q_k(D_xf^{-1}(\mathbf{s}_{fx}))-Q_k(D_xf^{-1}(\mathbf{s}_{fz}))} &<& 
\eps',
\text{ and } \\
\abs{Q_k(D_xf^{-1}(\mathbf{s}_{fz}))-Q_k(D_zf^{-1}(\mathbf{s}_{fz}))} &<& 
\eps'.
\end{eqnarray*}
Thus~$(2)$ follows from $Q_j(\mathbf{s}_{fz}) - 
Q_k(D_zf^{-1}(\mathbf{s}_{fz})) \geq d_3 - 
3\eps ' > 0$.
\end{proof}

Now we use a partition of unity to smooth $\{ Q_k \}_{k=1}^N$ on the invariant line 
fields.

\begin{defn} \label{defn:Qz}
Let $\Gamma$ be box hyperbolic.  Let $\tau  >0$ be 
as in Lemma~\ref{lem:deltaQexp}.
Define a partition of unity on $\mathcal{B}$ by choosing 
continuous functions $\rho_k:\Ct \to [0,1]$ for  each box $B_k
\in \mathcal V$, such that $\text{supp}(\rho_k) \subset \snbd{B_k}{\tau }$
and $\sum_{k} \rho_k(z) = 1$, for any $z\in \mathcal B$.

Let $z\in J$.  Then we define $Q_z \colon E^u_z \cup E^s_z \to \RR$ by
\[
Q_z(\mathbf{w}) = \sum_{k} \rho_k (z) Q_k(\mathbf{w}).
\]
\end{defn}

Note that $Q_z(\mathbf{w})$ is a continuous function of $\mathbf{w}$ since 
$Q_k$ is continuous, and further a continuous 
function of $z$ due to the partition of unity.

\begin{prop} \label{prop:combineQz}
Let $\Gamma$ be box hyperbolic.
Let $\{ Q_z \}$ be given by Definition~\ref{defn:Qz}.
Then for any $z \in J$ we have:
\begin{enumerate}
\item if $\mathbf{w}\in E^u_z$, then
        $Q_{fz}\( D_zf(\mathbf{w}) \) > Q_z (\mathbf{w})$;
\item if $\mathbf{v}\in E^s_{fz}$, then
         $-{Q_z\(\left[D_zf\right]^{-1}(\mathbf{v})\)} \ > \
     -{Q_{fz}(\mathbf{v})}.$
\end{enumerate}
\end{prop}

\begin{proof}
Let $\mathbf{u}_z\in E^u_z$ be such that $\norm{\mathbf{u}_z}=1$.
If we set
\begin{eqnarray*}
q^u_{f,z} &=& \min \{ {Q_j(D_zf(\mathbf{u}_z))} 
\colon f(z)\in \text{ supp}(\rho_j)  \},
\text{ \ and \ } \\
q^u_z &=& \max \{ {Q_k(\mathbf{u}_z)} \colon z\in \text{ supp}(\rho_k)\}, 
\end{eqnarray*}
then by Lemma~\ref{lem:deltaQexp} we know that $q^u_{f,z} > q^u_{z}$.  
Thus we need only 
use that the partition functions sum to one to get
\begin{eqnarray*}
Q_{fz}({D_zf(\mathbf{u}_z)}) 
& = & \sum_{j} \rho_j (f(z)) 
	{Q_j(D_zf(\mathbf{u}_z))}
\geq \sum_{j} \rho_j (f(z)) q^u_{f,z} = q^u_{f,z} \\
& > & q^u_{z} = \sum_{k} \rho_k(z) q^u_{z} \geq 
 \sum_{k} \rho_k(z) {Q_k(\mathbf{u}_z)}
= Q_z(\mathbf{u}_z).
\end{eqnarray*}
Hence $(1)$ follows since $Df$ is linear, and for any $\alpha \in \CC$,
$Q(\alpha \mathbf{w}) = \abs{\alpha}^2 Q(\mathbf{w})$.

Establishing $(2)$ is analogous.
Let $\mathbf{s}_{fz}\in E^s_{fz}$ be such that $\norm{\mathbf{s}_{fz}}=1$.
If we set
\begin{eqnarray*}
-q^s_{z} &=& \min \{ -{Q_k(\left[D_zf\right]^{-1}(\mathbf{s}_{fz}))} 
\colon z\in \text{ supp}(\rho_k)\}, 
\text{ \ and \ } \\
-q^s_{f,z} &=& \max \{ -{Q_j(\mathbf{s}_{fz})} 
\colon f(z)\in \text{ supp}(\rho_j)  \},
\end{eqnarray*}
then by Lemma~\ref{lem:deltaQexp} we know that $-q^s_{z} > -q^s_{f,z}$.  
Thus we need only 
use that the partition functions sum to $1$ to get
\begin{eqnarray*}
-Q_{z}({\left[D_zf\right]^{-1}(\mathbf{s}_{fz})}) 
 =  -\sum_{k} \rho_k (z) 
	{Q_k(\left[D_zf\right]^{-1}(\mathbf{s}_{fz}))}
\geq -\sum_{k} \rho_k (z) q^s_{z} = -q^s_{z} \\
 > -q^s_{f,z} = -\sum_{j} \rho_j(f(z)) q^s_{f,z} \geq 
 -\sum_{j} \rho_j(f(z)) {Q_j(\mathbf{s}_{fz})}
= -Q_{fz}(\mathbf{v}).
\end{eqnarray*}
\end{proof}

\begin{defn} \label{defn:pnorm}
Let $\Gamma$ be box hyperbolic.
Let $z\in J$. We define the norm $\pnorm{\cdot}{z}$ on $T_z\Ct$ using $Q_z$ 
and the spaces $E^u_z, E^s_z$ as a 
basis, \textit{i.e.}, for $\mathbf{w} \in T_z\Ct$,
\[
\pnorm{\mathbf{w}}{z} = \max \( 
\abs{Q_z(P^{u_z}_{s_z}(\mathbf{w}))}^{1/2}, 
\abs{Q_z(P^{s_z}_{u_z}(\mathbf{w}))}^{1/2} \),
\]
where $P^{a_z}_{b_z}$ denotes the projection onto 
$E^a_z$ with $E^b_z$ as its Null space.
\end{defn}

This is a continuous norm for $z \in J$.
Robinson (\cite{Rob}) notes in his construction of an adapted metric for hyperbolic sets that the maximum of two norms on subspaces defines a norm, which is very similar to the above.

Finally, we establish Theorem~\ref{thm:truehyptwo}, by 
showing that for a box hyperbolic $\Gamma$, we have $f$ is hyperbolic on $J$ with
respect to the norm $\pnorm{\cdot}{z}$ on $T_{z}\Ct$, for $z \in J$:

\begin{proof}
[Proof of Theorem~\ref{thm:truehyptwo}]
Suppose $\Gamma$ is box hyperbolic. We want to show $f$ is
hyperbolic over $J$, as in Definition~\ref{defn:hyp}, \textit{i.e.}, there is
a constant $\lambda > 1$, and for each $z$ in $J$ there is
a continuous splitting of the tangent bundle $T_{z}\Ct = E^s_z
\oplus E^u_z$, and a
continuous norm $\pnorm{\cdot}{z}$ such that:
\begin{enumerate}
\item $f$ preserves the splitting, \textit{i.e.}, $D_zf(E^s_z) = 
E^s_{fz}$, 
and
$D_zf(E^u_z) = E^u_{fz}$, and
\item $Df \ (Df^{-1})$ expands on $E^u_z \ (E^s_z)$ uniformly, 
\textit{i.e.},
 \begin{enumerate}
  \item if $\mathbf{w} \in E^u_z$ then
$\pnorm{D_zf(\mathbf{w})}{fz} \geq \lambda \pnorm{\mathbf{w}}{z}$, and
  \item if $\mathbf{w} \in E^s_z$ then
$\pnorm{D_zf^{-1}(\mathbf{w})}{f^{-1}(z)} \geq \lambda 
\pnorm{\mathbf{w}}{z}$.
  \end{enumerate}
\end{enumerate}

We have $(1)$ by Lemma~\ref{lem:splitting}.
Let $z\in J$ and $\pnorm{\cdot}{z}$ be given by Definition~\ref{defn:pnorm}.
 We show that $(2)$ follows easily from 
Proposition~\ref{prop:combineQz}.

First suppose $\mathbf{w} \subset E^u_z$.  Then $D_zf(\mathbf{w}) \subset
E^u_{fz}$. Hence $\pnorm{\mathbf{w}}{z}^2 = Q_z(\mathbf{w})$ and
$\pnorm{D_zf(\mathbf{w})}{fz}^2 = Q_{fz}(D_zf(\mathbf{w})).$  Thus 
Condition~$1.$ of Proposition~\ref{prop:combineQz} implies 
that $\pnorm{D_zf(\mathbf{w})}{fz} > \pnorm{\mathbf{w}}{z}$.

Now consider $\mathbf{w} \subset E^s_z$.  Then $
D_zf(\mathbf{w}) \subset E^s_{fz}$. 
Hence $\pnorm{\mathbf{w}}{z}^2 =  -Q_z(\mathbf{w})$ and
$\pnorm{D_zf(\mathbf{w})}{fz}^2 = -Q_{fz}(D_zf(\mathbf{w})).$ 
Then Condition~$(2)$ of Proposition~\ref{prop:combineQz} applied to 
$\mathbf{v} = \( D_zf^{-1}(\mathbf{w})\)$ implies 
$\pnorm{D_zf^{-1}(\mathbf{w})}{f^{-1}(z)} > \pnorm{\mathbf{w}}{z}$.

Finally, by compactness of $J$ the strict inequalities imply the 
existence of some constant $\lambda >1$, proving $(2)$.
\end{proof}

\subsection{Using linear algebra to characterize box hyperbolicity}

First recall that a Hermitian form $Q \colon \Ct \to \RR$ is associated to a Hermitian matrix $A$, so that $Q(\mathbf{w}) = \mathbf{w}^* A \mathbf{w}$.
Note that if $(k,j)\in \mathcal{E}$ is any edge in the graph $\Gamma$, then for any $z \in B_k$, $(Q_j \circ D_zf)$ is
also a Hermitian form, given by
\begin{equation} \label{eqn:QjDHmatrix}
Q_j (D_{z}f(\mathbf{w})) 
= \mathbf{w^*} ((D_{z}f)^* A_j (D_{z}f)) \mathbf{w}. 
\end{equation}

\begin{prop} \label{prop:conepresexp}
Suppose $\{Q_k\}$ are Hermitian forms with
${C}^u_k = $ $\{ \mathbf{w} \colon $ $ Q_k(\mathbf{w}) \geq 0 \}$
and  ${C}^s_k = \{ \mathbf{v} \colon Q_k(\mathbf{v}) < 0 \}$, for each box $B_k$ in $\mathcal{V}$.
Then $\Gamma$ is box hyperbolic (using $\{Q_k\}$) 
iff 
for every $B_k \in \mathcal{V}$, every $z \in B_k$, and
every edge $(k,j) \in \mathcal{E}$, we have $((Q_j\circ D_zf) - Q_k)$ is positive definite. 
\end{prop}

\begin{proof}
($\mathbf{\Leftarrow}$) We begin with the reverse implication.
Let $z\in B_k$ and $B_j$ be a box such that $(k,j)\in \mathcal{E}$.  Then
 $Q_j(D_{z}f)(\mathbf{w}) > Q_k(\mathbf{w})$, for
all ${z}\in B_k$ and all $\mathbf{w}\in\Ct$.

First consider the unstable cones. Suppose $\mathbf{w} \in 
{C}^u_k$, so by definition $0 < Q_k(\mathbf{w})$. But then by 
hypothesis, we get
\[
0 < Q_k(\mathbf{w}) < Q_j(D_{z}f(\mathbf{w})).  
\]
Thus $D_{z}f(\mathbf{w}) \in C^u_j$, so
the unstable cones are preserved by 
$D_{z}f$, and we have established Condition~1 of 
box hyperbolicity.

Next consider the stable cones. First, we show that stable cone
preservation follows from unstable cone preservation, since they are
complementary.  Indeed, above we showed that $Df$ preserves the unstable
cones, \textit{i.e.}, $D_zf({C}^u_k) \subset {C}^u_j$.  Hence, ${C}^u_k
\subset [D_{z}f]^{-1} ({C}^u_j)$. But by definition, $C^s = \Ct \setminus C^u$. 
Thus ${C}^s_k \supset [D_{z}f]^{-1}(C^s_j)$ and so the
stable cones are preserved by $Df^{-1}$.

Now let $\mathbf{v} \in {C}^s_j$, so that
\[
0 < -Q_j(\mathbf{v}) = -(Q_j \circ D_zf) 
\([D_{z}f]^{-1}(\mathbf{v})\).
\] 
Then since we have stable cone preservation under $Df^{-1}$, we also know 
that 
\[0 < -Q_k( [D_{z}f]^{-1}(\mathbf{v})).
\]

Combining this with the negative of the hypothesis establishes
Condition~2 of box hyperbolicity, \textit{i.e.},
\[
-{Q_k( [D_{z}f]^{-1}(\mathbf{v}))} > -{Q_j(\mathbf{v})}.
\]

\bigskip

$\mathbf{(\Rightarrow)}$ Now we prove the forward implication. 
Suppose $\Gamma$ is box hyperbolic, 
\textit{i.e.}, we 
have Conditions~1 and~2 of Definition~\ref{defn:boxhyp}. Let $(k,j) \in \mathcal{E}$, and $z \in B_k$.
We consider $\mathbf{w}$ in each of three regions to show 
$((Q_j\circ D_zf) - Q_k)$ is positive definite.

\textbf{Case 1:} Suppose $\mathbf{w} \in {C}^u_k$.
Then by definition $0 < Q_k(\mathbf{w})$.

Since box hyperbolicity implies 
the unstable cones are preserved by $Df$, we have that 
$D_{z}f(\mathbf{w}) \in {C}^u_j$, so 
$
0 < Q_j(D_{z}f(\mathbf{w})).
$

Then Condition~1 of box hyperbolicity gives us
\[
{Q_j(D_{z}f(\mathbf{w}))} \ > \ {Q_k(\mathbf{w})},
\]
hence $((Q_j\circ D_zf) - Q_k)$ is positive on $\mathcal{C}^u$.

\textbf{Case 2:} Suppose $\mathbf{w} \in [D_{z}f]^{-1}(\mathcal 
C^s_j)$, \textit{i.e.}, $\mathbf{v} = D_{z}f(\mathbf{w}) \in 
\mathcal 
C^s_j$. Then by definition $Q_j(D_{z}f(\mathbf{w})) <0.$

Now by stable cone preservation, we know $\mathbf{w} \in \mathcal 
C^s_k$, hence $Q_k(\mathbf{w}) <0$.  

Condition~2 of box hyperbolicity says that 
\[
-{Q_j(\mathbf{v})} \ < \ -{Q_k\([D_{z}f]^{-1}(\mathbf{v})\)}
\]
for all vectors $v\in C^s_j$, hence it applies to 
$\mathbf{v}=D_{z}f(\mathbf{w})$.  Thus we get
\[
-{Q_j\(D_{z}f(\mathbf{w})\)} \ < \ -{Q_k(\mathbf{w})},
\]
and negating yields 
\[
Q_j(D_zf(\mathbf{w})) > Q_k(\mathbf{w}),
\]
so $((Q_j\circ D_zf) - Q_k)$ is positive on $ [D_{z}f]^{-1}(\mathcal{C}^s)$.

\textbf{Case 3:} For the remaining $\mathbf{w}$, we have  $\mathbf{w} 
\notin {C}^u_k$ and $\mathbf{w} \notin 
[D_{z}f]^{-1}(C^s_j)$.  Then $Q_k(\mathbf{w})<0$ and 
$Q_j(D_{z}f(\mathbf{w})) \geq 0$. Hence,
\[
Q_j(D_{z}f(\mathbf{w})) \geq 0 > Q_k(\mathbf{w}).
\]
Thus we easily get 
$ Q_j(D_zf(\mathbf{w})) - Q_k(\mathbf{w}) > 0$.
\end{proof}

\section{Verifying box hyperbolicity: the \hyptwoalg }
\label{sec:hyptwoalg}

In this section, we explain in detail the {\hyptwoalg } for testing box hyperbolicity 
of a box chain model $\Gamma$ of a {\Henon } mapping, $H$, by attempting to construct a cone field and norm for which the map is hyperbolic.


But first, before we can test box hyperbolicity, we must start with a  $\Gamma$ which seems to model $H$ reasonably well.  Thus we now summarize how we use the {\boxchcn } of \cite{SLHone} to obtain separate strongly connected graphs modeling $J$ and any other invariant sets of recurrent dynamics, for example, sink cycles (attracting periodic orbits).  Recall from Section~\ref{sec:hypintro} that if $H$ is indeed hyperbolic, then the only recurrent dynamics are $J$ and a finite number of attracting periodic orbits. The construction is an iterative process.  We begin by defining a large box $\mathcal{B}_0 = [-R, R]^4$ in $\RR^4 = \Ct$, such that all the recurrent dynamics of the map is contained in $\mathcal{B}_0$ (in \cite{SLHone} we give a simple formula for $R$ in terms of the parameters $a,c$).  Then for some $n>1$, we place a $2^n \times 2^n \times 2^n \times 2^n$ grid of boxes on $\mathcal{B}_0$.  The construction then builds strongly connected graphs $\Gamma_1, \Gamma_2, \ldots, \Gamma_N$, each consisting of a subcollection of these grid boxes, and such that $J$ is covered by the boxes of $\Gamma_1$, and each sink cycle is covered by the boxes of some $\Gamma_k$.   
Then $\Gamma=\Gamma_1$ is a {\boxchmod } of $J$, as in Definition~\ref{defn:boxmodel}.  Using smaller boxes in the construction produces a more accurate {\boxchmod}.  

To prove hyperbolicity, we need each sink cycle in a different model from $J$, \textit{i.e.}, in some $\Gamma_k$ for $k\neq 1$.  If there seems to be a sink cycle together with $J$ in $\Gamma_1$, then we subdivide and repeat the above process. That is, place a grid of boxes inside of each box of $\Gamma_1$, and use these smaller boxes to obtain a refinement, $\Gamma_{1,1}, \Gamma_{1,2}, \ldots, \Gamma_{1,M}$, such that $\Gamma=\Gamma_{1,1}$ contains $J$.  If in this refinement, the sink cycle is in some $\Gamma_{k,1}$ for $k > 1$, then we can stop and study $\Gamma$.  Otherwise, repeat the subdivision process, until computational resources are exhausted, or a $\Gamma$ containing only $J$ is produced.

\begin{center}
\begin{figure}
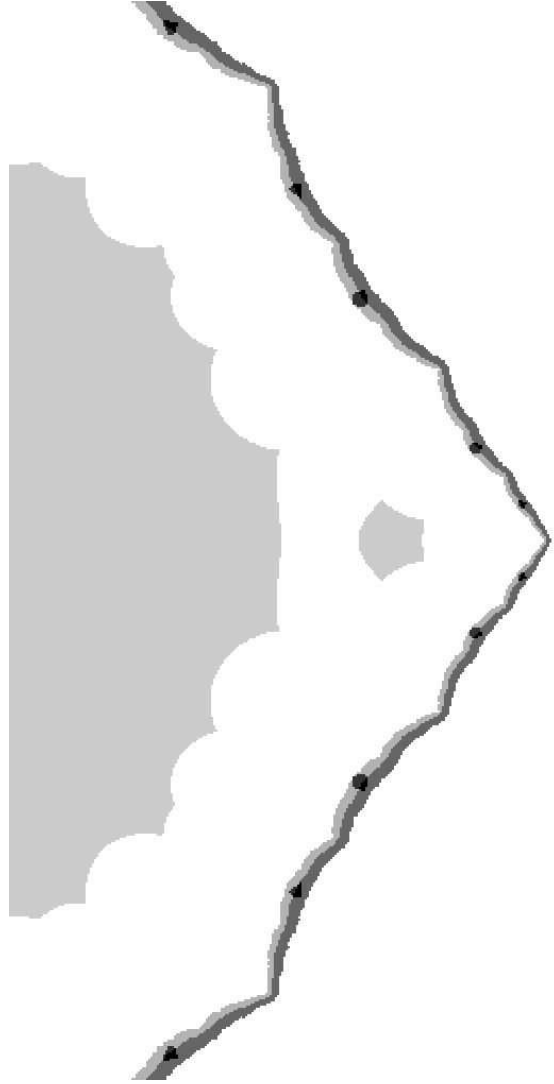

\drawfighyptwoperoneA
\caption[A box model of $J$ and the fixed sink (the only recurrent dynamics) for $H_{a,c}$, with $c=-.3, a=.1$]
{Shown above is the unstable manifold parameterization for a saddle fixed point of $H_{a,c},$ with $c=-.3, a=.1$, intersected with three box models $\Gamma_1, \Gamma_2, \Gamma_3$, containing pseudo-recurrent dynamics of $H$.
The two islands in the large interior on the left are in $\Gamma_2$, the model of the fixed sink. 
The right hand band is $\Gamma_1$,  containing $J$. The lighter band on $\Gamma_1$ is approximately contained in $K^+$.  The darkest spots overlapping $\Gamma_1$ show $\Gamma_3$, which contains pseudo-recurrence but no true recurrence.  Thus $\Gamma_3$ would be eliminated for smaller box size. Here boxes are of side length $(2\times1.43)/2^7=0.022$, and $\Gamma_1$ contains $32{,}000$ boxes.
This $\Gamma_1$ is box hyperbolic.} 
\label{fig:perone}
\end{figure}
\end{center}

We can check our accuracy at each level in the iterative process by
producing pictures of the current $\Gamma$'s boxes intersected with an unstable
manifold of a saddle periodic point. 
As discussed in
Section~\ref{sec:henwuppic}, we can  parametrize an
unstable manifold by a plane, then to determine the coloration of a
pixel, we check whether the pixel intersects some boxes of~$\Gamma$.  
Since the picture is a parametrization of a 
manifold which does not line up with the axes in $\Ct$, a pixel may hit
more than one box, and in more than one strongly connected component $\Gamma_k$. The user may also 
decide
to lighten the pixels which are heuristically found to be in
$K^+$, to check visually how close the model is to $J$.  
For example, for the {\Henon } mapping $H_{a,c}$, with $(c,a) = (-.3, .1)$, Figure~\ref{fig:perone} shows a parameterized unstable manifold intersected with the boxes in models $\Gamma_1, \Gamma_2, \Gamma_3$, with $\Gamma_1$ containing $J$, $\Gamma_2$ containing the fixed sink,  and $\Gamma_3$ containing pseudo-recurrence but no true recurrence  (thus $\Gamma_3$ would be eliminated for smaller box size). In this figure each $\Gamma_k$ is shaded differently, and pixels heuristically found to be in $K^{+}$ are lightened.



Now assume we have obtained a model $\Gamma$ which appears to contain $J$ but no other recurrent dynamics.  We outline below the {\hyptwoalg } for verifying box hyperbolicity on $\Gamma$, then describe each step in detail.

\begin{alg}[\hyptwoalg]  
\label{alg:hyptwoalg} 
\begin{enumerate}
\item Define an approximately invariant splitting; specifically, a
pair of ``unstable'' and ``stable'' vectors, $\mathbf{u}_k, \mathbf{s}_k$,
for each box, $B_k \in \mathcal V(\Gamma)$.  
\item Build an
``unstable'' metric which is (approximately) uniformly expanded by $DH$ on the set of
unstable directions, $\{\CC \mathbf{u}_k\}$, and a ``stable'' metric which
is (approximately) uniformly expanded by $DH^{-1}$ on the set of stable directions, $\{\CC
\mathbf{s}_k\}$. 
\item Use the directions and their metrics
to define the cones as Hermitian forms.  
\item Finally, check
whether $H$ preserves the cone field, and whether with
respect to the Hermitian forms, $DH \ (DH^{-1})$ is expanding on the
unstable (stable) cones. 
\end{enumerate}
\end{alg}

In this algorithm, we will need to move around on the graph $\Gamma$ using the following:

\begin{defn} \label{defn:arbor}
A directed graph $\Gamma_0$ is an {\em arborescence} if there is 
a root vertex $v_0$ so that for any other vertex $v$, 
there is a unique simple path from $v_0$ to $v$.  Such a graph is a 
tree, and must have exactly one incoming edge for each vertex $v \neq 
v_0$.  

If $\Gamma$ is strongly connected, then for each vertex 
$v_0$ in $\Gamma$, there is a minimum spanning tree $\Gamma_0$ with 
root vertex $v_0$ which is an arborescence 
(simply perform a depth first or breadth first search from $v_0$).  We 
call such a $\Gamma_0$ a {\em spanning arborescence}.
(See \cite{TCR} for a discussion of implementation of basic graph theory).
\end{defn}


\begin{step}[Setting stable and unstable directions]
\label{alg:USdirs} 

Recall from Section~\ref{sec:henwuppic}, that when $\abs{a} \neq 1$, except on the curve of equation $4c = (1+a)^2$, the map $H_{a,c}$ has at least one saddle fixed point, $p$.   If $H$ has no saddle fixed point, it should be possible to instead use a saddle periodic point of period greater than one, but this was not necessary for any of the maps we were interested in studying.
The {\Henon } mapping $(x,y) \to (x^2+c-ay,x)$ has two fixed
points.  Note first that fixed points of~$H$ must be on the diagonal,
\textit{i.e.}, $x=y$.  Then the two fixed points are:
\[
x_{\pm} = \frac{1}{2} \( (1+a) \pm \sqrt{(1+a)^2+4c} \).
\]
By substituting the IA Hull of $a$ and $c$ into the above formula, and performing operations in IA, we compute $x_{\pm}$ as complex intervals (interval vectors in $\Rt$) containing the actual fixed points. Next, the eigenvalues of~$D_{x}H$ for the fixed point $(x, x)$ 
are:
\[
\lambda, \mu = x \pm \sqrt{x^2 - a}.
\]
Thus we test to see whether each interval vector $x_{\pm}$ is a saddle, by computing two real intervals containing the moduli of the eigenvalues, and then testing whether one interval is entirely greater than one, while the other is entirely less than one.
So now assume $H$ has a saddle fixed point $p=(x,x)$.  Bedford and Smillie have shown $J$ contains all saddle periodic points, hence $p\in J$, thus there is a box, $B_0$, in $\Gamma$, containing $p$.  Let $\mathbf{u_0} \ ( \mathbf{s_0} )$ be the eigenvector of
unit length corresponding to the unstable (stable) eigenvalue.  These are natural unstable and stable directions in the box $B_0$.  Since we will just use the $\mathbf{u_k}$'s and $\mathbf{s_k}$'s as axes for cones, we need only know them approximately, so here interval arithmetic is not needed.

Now let $\Gamma_0$ be a spanning arborescence of $\Gamma$, with root vertex $B_0$.
To define unstable directions in each box $B_k$, sucessively push
$\mathbf{u_0}$ across the edges of  $\Gamma_0$ by $DH$.  To be precise, fix a point $z_k$ in $B_k$, say the center point of the box, and for each edge $(k,j) \in \mathcal{E}({\Gamma_0})$, starting with $k=0$, define
\[
\mathbf{u}_{j}  := \frac{D_{z_k}H (\mathbf{u}_k)}{\| D_{z_k}H (\mathbf{u}_k) \| }.
\]  

As noted above, we need only an approximation to the $\mathbf{u_k}$'s, so interval arithmetic is not used here.  In fact, were we to attempt to use IA, we may encounter a computational difficulty.  This is due to the recursive definition of the $\mathbf{u_k}$'s.  Suppose we desired to compute ``more'' of the potential unstable directions in each box, by starting with an interval vector $\mathbf{U_0}$ in $\Ct$ guaranteed to contain the unstable eigenvector of $D_pH$, and then pushing this interval vector across $\Gamma_0$  to define interval vectors $\mathbf{U_k}$ in each box $B_k$.   More precisely, a box $B$ is an interval vector in $\Ct$, so $B = (X, Y)$, with $X, Y$ complex intervals (interval vectors in $\Rt  =\CC)$. Then 
\[
D_pH = \left[\begin{array}{rr}
2p &-a \\
1 &0 
\end{array}\right]
 \text{ and } D_{B} H = \left[ \begin{array}{rr}
 2X & -\text{Hull}(a) \\
 \text{Hull}(1) & \text{Hull}(0) \end{array} \right],
\] 
so $D_B H$ has entries complex intervals.  Then $D_{B_k} H (\mathbf{U_k})$ is an interval vector.  Suppose $(1, 2, \ldots, n)$ is a path in $\Gamma_0$ (this should be $(k_1, \ldots, k_n)$, but we wish to avoid too many nested subscripts). Then 
$
\mathbf{U_{n}} = D_{B_n} H\circ \cdots \circ D_{B_1} H(\mathbf{U_1}).
$ 
However, since the $B_k$'s are near $J$, due to the dynamical properties of $H$ on $J$, this type of iteration will (in experiment, very quickly) lead to intervals $\mathbf{U_n}$ to huge to be useful.

Define the stable directions similarly, keeping in mind that 
stable cones should be expanded and preserved by $DH^{-1}$.  
The transpose of a graph $\Gamma$, $\Gamma^T$, is the graph formed by 
reversing the edge directions of $\Gamma$.
Thus to define $\mathbf{s}_k$, we use a spanning tree $\overline{\Gamma_0}$ 
of~$\Gamma^{T}$ with $B_0$ as root vertex (note that $\overline{\Gamma_0} \neq \Gamma_0^T$), and push $\mathbf{s}_0$
across successive edges of $\overline{\Gamma_0}$ by $DH^{-1}$.  Specifically, if $(j,k) \in
\mathcal{E}({\overline{\Gamma_0}})$, then
\[
\mathbf{s}_{k}  =  \frac{D_{z_k}H^{-1} (\mathbf{s}_j)}{\| D_{z_k}H^{-1} (\mathbf{s}_j) \| },
\]
 where $z_k \in B_k$ as in the  definition of $\mathbf{u}_k$. 

\end{step}

\begin{prop} \label{prop:USdirInCone}
Let $\{\mathbf{u}_k\}, \{\mathbf{s}_k\}$ be unstable and stable directions 
for each box $B_k$ of $\Gamma$, defined as above. 
\begin{enumerate}
\item
Let ${C}^u_k$ be any box constant cone field preserved by $D_zH$,
for each $z\in B_k$. Then for each $B_k$, we must
have $\mathbf{u}_k \subset {C}^u_k$.
\item
Let ${C}^s_k$ be any box constant cone field preserved by 
$D_zH^{-1}$, for each $z\in B_k$. Then for each $B_k$, we must
have $\mathbf{s}_k \subset {C}^s_k$.
\end{enumerate}
\end{prop}


Note that for edges $(k,m)$ in $\Gamma$ but not in the spanning arborescence $\Gamma_0$, $DH$ does not 
map $\mathbf{u}_k$ into $\CC\mathbf{u}_m$.  It is helpful in establishing
invariance of the cone field if $DH(\mathbf{u}_k)$ is close to
$\CC\mathbf{u}_m$, and in addition, if $D_{z}H$ does not vary greatly as ${z}$ varies within
one box.

Thus before performing the next step, we take some measurements on the separation of the 
stable and unstable directions in each box, to get an idea of whether 
it might be possible to prove box hyperbolicity using these 
directions, and store (for later scrutiny) which boxes might be obstructions.
In order to measure the difference between directions, 
we view a direction in $\Ct$ as a complex line in $\Ct$, and thus use 
the spherical metric, $\sigma$.
Then, for each box $B_k$ in $\mathcal{V}(\Gamma)$, let
\begin{eqnarray*}
\text{Udiam}[k] & = &  \text{diam}_{\sigma} 
 	\{ D_{z_j}H(\mathbf{u_j}) \colon (j,k) \in \mathcal{E}(\Gamma) \}, 
   \text{ and} \\
\text{Sdiam}[k] & = &  \text{diam}_{\sigma} 
 	\{ [D_{z_k}H]^{-1}(\mathbf{s_j}) \colon (k,j) \in \mathcal{E}(\Gamma) \}.
\end{eqnarray*}  
We do not measure the
variation within one box, \textit{i.e.}, between $D_zH$ and $D_{z_k}H$ 
for different $z$ in box $B_k$,
since it seems that would be much smaller than among images
from different boxes.

Proposition~\ref{prop:USdirInCone} suggests that
a clear separation between Udiam$[k]$ and Sdiam$[k]$
is needed in order for a computer program to verify cone preservation 
under $DH$, thus we will not confidently progress to the next step unless we 
have
\[
\sigma(\mathbf{u}_k,\mathbf{s}_k) - (\text{Udiam}[k]+\text{Sdiam}[k])>0
\] 
in each box $B_k$.  Finally, we record in which boxes the
 either Udiam or Sdiam is large, or $\sigma - (\text{Udiam} + \text{Sdiam})$ is negative.


\begin{step}[Building a metric on the directions] 
\label{alg:buildmetdir}
Consider the unstable directions $\{\CC \mathbf{u}_k\}_{k=1}^N$.  As 
discussed above, $DH$ does not quite preserve these directions, so first we take
that into account.  Let $P^{\mathbf{u}}_{\mathbf{s}}$ be the projection
onto $\CC \mathbf{u}$ with $\CC \mathbf{s}$ as its Null space.  Given the
vectors $\mathbf{u} = (u_1, u_2)$ and $\mathbf{s} = (s_1,s_2)$ in $\Ct$,
this projection is:
\[
P^{\mathbf{u}}_{\mathbf{s}} \twovec{v_1}{v_2}
= \frac{v_1 s_2 - v_2 s_1}{u_1 s_2 - u_2 s_1} \twovec{u_1}{u_2}.
\]
Then for each edge $(k,j) \in \mathcal{E}(\Gamma)$, we have $P^{\mathbf{u}_j}_{\mathbf{s}_j} \circ D_{z_k}H$ maps $\CC
\mathbf{u}_k$ to $\CC \mathbf{u_j}$. If the unstable directions are a
good approximation to an invariant unstable line field, then
$P^{\mathbf{u}}_{\mathbf{s}} \circ DH$ is close to $DH$ on these unstable
directions.

In \cite{SLHtwo}, we describe a method for proving hyperbolicity of polynomial maps of one complex variable, by building a metric for which the map is expanding by some $L>1$ on a neighborhood of the Julia set.  The neighborhood of $J$ is a collection of boxes in $\CC$, and the metric in each box is defined by a constant times the $L^{\infty}$ metric in $\Rt=\CC$. Say $\varphi_k \norm{\cdot}$ is the metric on box $B_k$, then the constants $\varphi_k$ are called \textit{metric handicaps}. 

In two variables, we use this same algorithm twice: once for the unstable directions and once for the stable directions.  For example, for the unstable direction field, we will attempt to build a metric for which
$P^{\mathbf{u}}_{\mathbf{s}} \circ DH$ is box expansive by some amount
$L>1$.   This metric will be defined in each unstable direction $\CC \mathbf{u_k}$ by a constant times our $L^{\infty}$ metric in $\RR^{4} = \Ct$, say $\varphi^u_k \norm{\cdot}$.  Following the algorithm in \cite{SLHtwo}, we call the $\varphi^u_k$ \textit{unstable
metric handicaps}, and define $\varphi^u_0=1$ in box $B_0$, then want to build handicaps satisfying
\[
\varphi^u_j  {\norm{P^{\mathbf{u}_j}_{\mathbf{s}_j} \circ D_{z_k}H (\mathbf{v})}}
\geq  {L } \ \varphi^u_k \norm{\mathbf{v}} ,
\]
for each edge $(k ,j) \in \mathcal{E}(\Gamma)$, and each $\mathbf{v} \in \CC \mathbf{u_k}$.
But then $\mathbf{v} = \alpha \mathbf{u_k}$ for some $\alpha \in \CC$, so 
since $P^{\mathbf{u}}_{\mathbf{s}} \circ DH$ is linear and 
$\norm{\mathbf{u}_k}=1$, 
the above equation is equivalent to
\begin{equation} \label{eqn:umet}
\varphi^u_j \geq  \frac{L \ \varphi^u_k}{\norm{P^{\mathbf{u}_j}_{\mathbf{s}_j} \circ 
D_{z_k}H (\mathbf{u}_k)}}.
\end{equation}
If we set $\xi_{k,j} = L/ {\norm{P^{\mathbf{u}_j}_{\mathbf{s}_j} \circ D_{z_k}H (\mathbf{u}_k)}}$, then we can use Algorithm~4.4 of \cite{SLHtwo} to attempt to find metric handicaps satisfying $\varphi^u_j \geq \xi_{k,j} \varphi^u_k$, hence Equation~\ref{eqn:umet}, on each edge $(k,j) \in \Gamma$. 
 This is of course not always possible, but the intuition
is that it should be possible if the box model is sufficiently small in
the right places.  If it is not possible, then the user can try a smaller $L$, or start over with smaller boxes.

We use interval arithmetic in Equation~\ref{eqn:umet} as in \cite{SLHtwo}.  The $\varphi_k^u$'s will be chosen  in $\FF$ (machine knowable numbers, note $\varphi^u_0=1 \in \FF$).  Then replace each number in the left hand side with it's Hull (the smallest interval in $\IF$  containing it).  Then perform the operations in IA, and test whether the upper endpoint of the resulting interval is less than $\varphi^u_j$.

To attempt to define stable metric handicaps, use
the method analogous to that for the unstable metric handicaps.  That is,
try to find an $M>1$ and build handicaps $\{ \varphi^s_k\}$ so that
\begin{equation} \label{eqn:smet}
\varphi^s_j \leq \frac{\varphi^s_k M}{\norm{P^{\mathbf{s}_j}_{\mathbf{u}_j} \circ 
D_{z_k} H (\mathbf{s}_k)}},
\end{equation}
for each edge $(k ,j) \in \mathcal{E}(\Gamma)$. Then the stable directions are box-contracted by
$P^{\mathbf{s}}_{\mathbf{u}} \circ DH$.  Again, if this step fails to produce a contracted metric, then the user can try a larger $M$, or start over with smaller boxes.

As we will see below, the ratio of stable to unstable metric handicaps in each box, $\varphi^s_k$ to $\varphi^u_k$, determines the width of the cones. Hence it is necessary to find values for $L$ and $M$ which yield comparable metrics.  
Lyapunov exponents give
us some intuition.
For a polynomial automorphism of $\Ct$, there are 
two Lyapunov exponents, $\lambda^{\pm}$, which measure
expansion and contraction of tangent vectors.
According to \cite{BS6}, if $f$ is a polynomial diffeomorphism of $\Ct$ with $d=d(f)>1$,
then $\lambda^{+} \geq \log d$, $\lambda^{-} \leq -\log d$, and
\begin{equation} \label{eqn:lyaptwo}
\lambda^+ + \lambda^- = \log(\det Df).
\end{equation}  
Note that for {\Henon } mappings, ${\det(Df)} = a$. Thus
$\lambda^+ + \lambda^- = \log(a).$ Since $\lambda^{+} \geq \log d$, we
have the inequality:  $\lambda^- \leq \log(a) - \log d$. In the case
$\abs{a} < 1$, we have $\log(a) <0$, hence the inequality for $\lambda^-$
is stronger than the inequality for $\lambda^+$. Thus in general we expect
stronger contraction than expansion of tangent vectors under {\Henon }
mappings.  
Equation~\ref{eqn:lyaptwo} implies that a good rule of thumb for 
choosing $L$ and $M$ is $L M = \abs{a}$.   In practice, for any $\Gamma$ we tested, the algorithm for setting the stable metric handicaps given any $M$ always completed in much less time than the algorithm for setting unstable metric handicaps given an $L$ (perhaps due to the strong contraction).
Thus after experimenting with various options, we have adopted the following strategy. First find the smallest working $M$ using simple bisection ($0<M<1$), then test expansion on the unstable directions using a value of $L$ near $\abs{a}/M$.  

In experiment we have observed that finding a working $L>1$ and $M<1$ is almost always possible when a box model $\Gamma$ has been found for $J$ which does not contain any sinks of $H$.  Rather, the difficult step is the next one: checking whether the cone field defined by these metrics is preserved and expanded by $DH$.
\end{step}


\begin{step}[Defining a cone field]
\label{alg:defcones}
If Step~\ref{alg:buildmetdir} successfully constructed expanded and contracted metrics on the unstable and stable directions, respectively, then the metrics and directions always define cones in each box, as follows.

In each box, $B_k$, define the unstable cone, $C^u_k$, so that a
vector $\mathbf{w}$ is in the unstable cone if it is closer to
$\mathbf{u}_k$ than $\mathbf{s}_k$, relative to the unstable and stable
metrics.  That is, $\mathbf{w} \in C^u_k$ iff 
\[
\varphi^u_k
\norm{P^{\mathbf{u}_k}_{\mathbf{s}_k}(\mathbf{w})} \geq \varphi^s_k
\norm{P^{\mathbf{s}_k}_{\mathbf{u}_k}(\mathbf{w})}.
\] 
Then the stable cones are just
the complements, $C^s_k := \Ct \setminus C^u_k$.

We define the Hermitian
form $Q_k \colon \Ct \to \RR$, by
\[
Q_k(\mathbf{w}) = \(\varphi^u_k 
\norm{P^{\mathbf{u}_k}_{\mathbf{s}_k}(\mathbf{w})}\)^2 
- \(\varphi^s_k \norm{P^{\mathbf{s}_k}_{\mathbf{u}_k}(\mathbf{w})}\)^2.
\]
Thus the unstable cone, $C^u_k$, is simply the set of vectors 
for which $Q_k$ is nonnegative, and the stable cone, $C^s_k$, is the
set of vectors for which the form is negative.

We can construct a Hermitian matrix, $A_k$, which encodes the information
of~$Q_k$, following standard linear algebra as in \cite{HK}.  A Hermitian
form $Q$ defines a sesquilinear form $g \colon \Ct \times \Ct \to
\RR$, such that $g(\mathbf{w},\mathbf{w}) = Q(\mathbf{w})$, where we can
recover $g$ using:
\[
g(\mathbf{v},\mathbf{w})  
= \frac{1}{4} Q( \mathbf{v} + \mathbf{w} )
-  \frac{1}{4}  Q( \mathbf{v} - \mathbf{w} ) 
+ \frac{i}{4} Q( \mathbf{v} + i \mathbf{w} )
-  \frac{i}{4} Q( \mathbf{v} - i \mathbf{w} ).
\]

A sesquilinear form $g$ can be represented by a matrix $A$ so that
$g(\mathbf{v},\mathbf{w}) = \mathbf{w^*} A \mathbf{v}$, with $a_{m,n} =
g(\mathbf{e}_n,\mathbf{e}_m)$ for an ordered basis $\{ \mathbf{e}_1,
\mathbf{e}_2 \}$, like $\{ (1,0), (0,1) \}$. Now $g$ Hermitian implies
that $A$ is also Hermitian, and the range of~$Q$ is $\RR$.  Thus, 
$Q(\mathbf{w})
= \mathbf{w^*} A \mathbf{w}$, where
$a_{m,n} = \frac{1}{4} Q( \mathbf{e}_n + \mathbf{e}_m )
-  \frac{1}{4}  Q( \mathbf{e}_n - \mathbf{e}_m ) 
+ \frac{i}{4} Q( \mathbf{e}_n + i \mathbf{e}_m )
-  \frac{i}{4} Q( \mathbf{e}_n - i \mathbf{e}_m )$.

We easily calculate that for $\mathbf{u} = (u_1, u_2), \mathbf{s} = (s_1,
s_2)$, if we set
\[
b_{11} = (\varphi^u \abs{s_2} \norm{u})^2 - (\varphi^s \abs{u_2} \norm{s} )^2, \ \ 
b_{22}  =  (\varphi^u \abs{s_1} \norm{u})^2 - (\varphi^s \abs{u_1} \norm{s} )^2, 
\]
\begin{eqnarray*}
b_{12} & = \frac{1}{4} & \left[
(\varphi^u \norm{u})^2 (\abs{s_2-s_1}^2 - \abs{s_2+s_1}^2 
+ i \abs{i s_2-s_1}^2 - i \abs{i s_2+s_1}^2) \right. \\
&  & 
- \left. (\varphi^s \norm{s})^2 ( \abs{u_2-u_1}^2 - \abs{u_2+u_1}^2 
+ i \abs{i u_2-u_1}^2 - i \abs{i u_2+u_1}^2) \right] 
,
\end{eqnarray*}
and $b_{21} = \bar{b}_{12}$,
then $a_{m, n} = b_{m, n}/{\abs{u_1 s_2 - u_2 s_1}^2}$.

In implementation, we use the above formulas to calculate, for each box $B_k$ in $\Gamma$, an interval  valued matrix $A_k$, representing the form $Q_k$ defining the cone $C^u_k$.  Note that since $A_k$ is Hermitian, the main diagonal entries are real intervals, and the other entries are (complex conjugate) complex intervals.    As noted above, the vectors $\mathbf{u_k}$ and $\mathbf{s_k}$ and the handicaps $\varphi^u_k, \varphi^s_k$ are chosen to be machine knowable numbers. Thus before the formulas defining $A_k$ are computed, these terms are converted to their interval Hulls, then the arithmetic operations are performed in IA to obtain $A_k$ with interval entries (of length greater than zero).
\end{step}

\begin{rem} \label{rem:thincones}
Note that the ratio of the metric handicaps determines the angle width of
the cones.  Thus if $\varphi^u_k$ and $\varphi^s_k$ are several orders
of magnitude different then the cones will be very thin, even if the unstable
and stable directions are far apart, thus the cones will be difficult
for the computer to work with.  This is why it is necessary in 
Step~\ref{alg:buildmetdir} to find values of $L$ and $M$ which yield
a comparable pair of metric handicaps in each box.
\end{rem}


\begin{step}[Checking whether cones are preserved and expanded]
\label{alg:checkcones}
For the last step of testing box hyperbolicity, we need to test whether
$DH \ (DH^{-1}$) expands the unstable (stable) cones, with respect to $\{
Q_k \}$. For this step we simply use Proposition~\ref{prop:conepresexp}, in
which we showed that in order
 to get preservation and expansion of the unstable cones,
and contraction of the stable cones,
we need precisely
 that $((Q_j\circ D_zH) - Q_k)$ is positive 
definite for every edge $(k,j) \in
\mathcal{E}(\Gamma)$, and every $z \in B_k$.  Thus in this step, we simply compute
this form defined for each edge in the graph, and test whether it is positive definite. 

In Step~\ref{alg:defcones}, for each box $B_k$ in $\Gamma$, we computed an interval matrix $A_k$ representing $Q_k$, in that $Q_k(\mathbf{w}) = \mathbf{w}^* A_k \mathbf{w}$. Thus by
Equation~\ref{eqn:QjDHmatrix}, the interval matrix representing 
$
(Q_j\circ D_zH) - Q_k
$
is 
\[
T_{k,j} = ((D_{B_k}H)^* \ A_j \ D_{B_k} H) - A_k.
\] 
Using the formula for the interval matrix $D_{B_k} H$ from Step~\ref{alg:USdirs}, it is straightforward to compute the interval matrix $T_{k,j}$.  Now we need to check whether this matrix is positive definite.  But since $T$ is Hermitian, the trace and determinant are real, and $T$ is positive definite iff its trace and determinant are positive (see \cite{HK}).  The trace and determinant of $T$ are real intervals, so we simply compute them with IA and check whether their lower endpoints are positive.

If the above test succeeds (positive definite for each edge), then the model $\Gamma$ is box hyperbolic, hence $H$ is hyperbolic.  If not,  then 
we may record boxes which are obstructions, that is, boxes $B_k, B_j$ for which $T_{k,j}$ fails to be positive definite. 
\end{step}

This is the stopping point of the {\hyptwoalg } for testing box hyperbolicity for a given $\Gamma$.  If box hyperbolicity fails, the user may refine $\Gamma$ by choosing to subdivide either all the boxes, or some subset of 
the boxes which seem to be obstructing the hyperbolicity testing (for 
example, boxes marked in Steps~\ref{alg:USdirs} or~\ref{alg:checkcones} above), then test the new $\Gamma$ with the {\hyptwoalg}.  
Two of the {\Henon } mappings of Theorem~\ref{thm:examples}, $H_{a,c}$ with $(c,a) = (0, -.22)$ and $(1.5, .5)$, were proven hyperbolic for a model constructed 
by straight subdivision to a certain box size, then by subdividing twice only 
boxes marked in Step~\ref{alg:checkcones} (see next section).


\section{Results of running {\Hypatia } on {\Henon } mappings}
\label{sec:examples}

Running {\Hypatia } for a {\Henon } mapping is not quite as 
simple as inputting the parameters $a,c$ and awaiting the results, since the user 
must make decisions as to how to build the best $\Gamma$ for testing with 
the {\hyptwoalg}.  In this section, we describe the 
specific process we followed and results obtained for the mappings of Theorem~\ref{thm:examples}.


%
%
%

Theorem~\ref{thm:examples} follows from Theorem~\ref{thm:truehyptwo}, that box hyperbolicity of some $\Gamma$, a {\boxchmod } of $J$,  implies hyperbolicity of $H$ on $J$, and from the fact that for each mapping mentioned in the theorem, using our computer program {\Hypatia}, we constructed a {\boxchmod } of $J$ and verified box hyperbolicity with the {\hyptwoalg}. Below, we discuss the process for each of these mappings in increasing order of the computational difficulty of proving hyperbolicity. 

%
The quickest map to be proven hyperbolic was $H_{a,c},$ with $c=-.3, a=.1$. 
We 
simply used a {\boxchmod } of $J$ with boxes selected from an
evenly subdivided $(2^7)^4$ grid on $\mathcal{B_0}=[-1.43, 1.43]^4$. Figure~\ref{fig:perone} shows the
box hyperbolic $\Gamma$. This is a map seemingly in the main cardioid, with recurrent dynamics $J$ and a fixed sink.

\begin{figure}
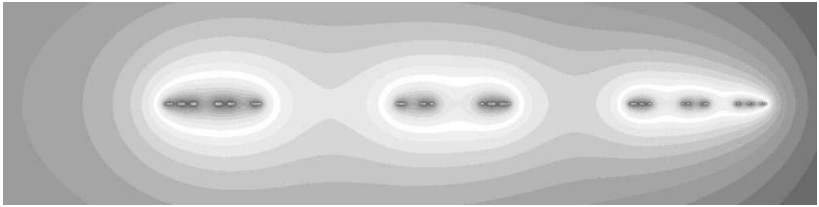

\begin{center}
\drawfighyptwocantorone
\caption{A FractalAsm picture of $W^u_p \cap K^+$ for $H_{a,c},$ with 
$c=-3, a=-.25$.  This set appears to be a real Cantor set.}
\label{fig:cantorone}
\end{center}
\end{figure}

For $H_{a,c},$ with 
$c=-3, a=-.25$, we also proved hyperbolic relatively quickly. The {\boxchmod } $\Gamma$ of $J$ from an 
evenly subdivided $(2^7)^4$ grid on $B_0=[-2.57, 2.57]^4$ is box hyperbolic.  This mapping appears to be a real horseshoe (\textit{i.e.}, a horseshoe contained in $\Rt$).
 Figure~\ref{fig:cantorone} is a
FractalAsm picture of the Julia set.  This kind of
picture is really the most useful for a Cantor set.

We proved the map $H_{a,c},$ with $c=0, a=-.22$, is hyperbolic by starting with a model of $J$ from an 
evenly subdivided $(2^7)^4$ grid on $\mathcal{B}_0=[-R,R]^4$, for $R=1.32$, but then additionally performing three 
hyperbolicity tests, and each time subdividing only boxes 
in which the cone check of Step~\ref{alg:checkcones} of the {\hyptwoalg } failed (to end with boxes of size ranging from 
$2R/2^7$ to $2R/2^{10}$).  This map also seems to be in the main cardiod.  The picture of the Julia set is similar to Figure~\ref{fig:perone}.

\begin{figure}
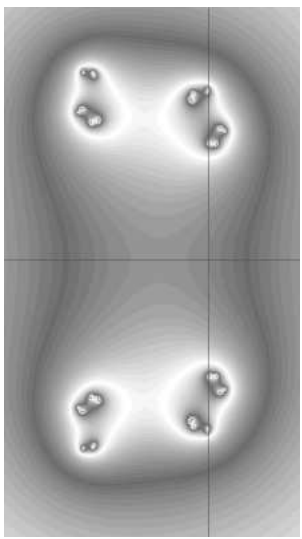

\begin{center}
\drawfighyptwocantortwo
\caption{A FractalAsm picture of $W^u_p \cap K^+$ for $H_{a,c},$ with 
$c=1.5, a=.5$.  This set appears to be a complex Cantor set.  The origin of the axes is the saddle fixed point $p$.}
\label{fig:cantortwo}
\end{center}
\end{figure}

Using nearly the same method as the previous mapping, we proved $H_{a,c},$ with $c=1.5, a=.5$, is hyperbolic.  Here, we started with the even 
$(2^8)^4$ grid on $\mathcal{B}_0=[-R,R]^4$, for $R=2.286$, then twice subdivided only boxes in which the cone check (Step~\ref{alg:checkcones} of the \hyptwoalg)
failed (yielding boxes of size $2R/2^8$ to $2R/2^{10}$).  The resulting {\boxchmod } is box hyperbolic.  FractalAsm pictures (see Figure~\ref{fig:cantortwo}) suggest this map is a complex horseshoe, with Julia set not contained in $\Rt$.

Table~\ref{table:examples} contains more data for all of the mappings discussed in this section.  In the table, $\Gamma$ denotes the {\boxchmod } of $J$, for the map $H_{a,c} (x,y) = (x^2+c-ay,x)$.  $R$ is the bound such that the boxes are contained in $\mathcal{B}_0 = [-R,R]^4 \subset \RR^4=\Ct$. The \textit{box grid depth} for a box is the number $n$ such that the box is of size $2R/2^n$.  If a {\boxchmod } contains boxes of multiple sizes, then multiple box grid depths are listed.

\begin{sidewaystable}
\caption{Data for the {\boxchmods } of Theorem~\ref{thm:examples}. [*Data not collected.]}
\label{table:examples}
\begin{center}
\begin{tabular} {| |    l  l  ||c|c|c |c ||}
\hline \hline

Figure & & \ref{fig:cantortwo} & \ref{fig:cantorone} & \ref{fig:perone}     &N/A 
\\ \hline

params. & $c$ & $1.5$  &      $-3$   &   $-.3$   &  $0$  
 \\ \cline{2-6}
	      & $a$  & $.5$    &     $-.25$  &   $.1$    &    $-.22$ 
\\ \hline
sink period &        &  N/A  &  N/A  &  $1$  & $1$ 
\\ \hline
$R$    &   &  $2.286$   &  $2.57$    &  $1.43$   & $1.32$ 
\\ \hline
box grid depth, $n$  & & $8$-$10$  & $7$  &  $7$  & $7$-$10$    
\\ \hline
 box size & &  $0.018, 0.009, 0.0045$   & $0.04$ & $0.022$   &$0.0052, 0.0026, 0.0013$ 
\\  \hline
$\Gamma$ size ($1000$s) &   boxes &  $53$   & $2.4$        & $32$  &  $78$  
\\ \cline{2-6}
           &  edges &  $2{,}950$   & $75$ & $1{,}250$ & $4{,}150$  
\\ \hline
$\sigma(\mathbf{u},\mathbf{s})$   (in $[0,\pi/2]$) &    min. &  *  &  $0.47$   &  $ 0.9$   & * 
 \\ \cline{2-6}
				    &    avg. & *   & $1.15$ &  $1.07$  & * 
\\ \hline
Udiam 	(in $[0,\pi/2]$)	   & max. &  * & $0.13$ & $0.015$ & * 
\\ \cline{2-6}
 & avg. &  * & $0.023$ & $0.03$ & * 
 \\ \hline
Sdiam 	(in $[0,\pi/2]$)	   & max. &  * & $0.11$ & $0.0008$ & *  
\\ \cline{2-6}
 & avg. &  * & $0.01$ & $0.0003$ & *  
 \\ \hline
$\sigma - (\text{Udiam}+\text{Sdiam}) $ & min.     
& * & $0.47$ & $.09$ & *  
\\ \cline{2-6}
 (in $[0,\pi/2]$)						& avg. 
&  * & $1.15$ & $1.07$ & *  
\\  \hline
$M$  (bisection) &  &  $0.21680$ & $0.111328$ & $0.0596$ & $0.140$ 
\\ \hline
$L$  $(\approx \abs{a}/M$) &  &  $2.180$ & $2$ & $1.6789$ & $1.480$  
\\ \hline
$\varphi^s_k $  & min. &  $0.0613$ & $0.046$ & $0.034$ & $0.0425$  
\\ \cline{2-6}
 & max. &  $1.079$ & $1$ & $1$ & $1$ 
 \\ \cline{2-6}
    & avg. &  $0.571$ & $0.096$ & $0.152$ & $0.263$ 
\\ \hline
$\varphi^u_k $  (in $(0,1]$) & min. &  $0.046$ & $0.033$ & $0.072$ & $0.0574$ 
\\ \cline{2-6}
        & avg. &  $0.268$ & $0.09$ & $0.158$ & $0.177$  
\\ \hline
proved  box-hyp? &   &  YES & YES & YES & YES  
\\ \hline
runtime (min.)   & &    $<10$    &  $<1 $ & $30$   & $45$    
\\  \hline
RAM (MB)          & & $>250$      & $20$ & $220$ & $1200$  
\\ \hline

\hline
\hline
\end{tabular}
\end{center}
\end{sidewaystable}

\bibliographystyle{plain}
\bibliography{lynch}

\def\cprime{$'$}
\begin{thebibliography}{10}

\bibitem{CUweb}
Dynamics at~Cornell.
\newblock [http://www.math.cornell.edu/\~{}dynamics].

\bibitem{BLS1}
E.~Bedford, M.~Lyubich, and J.~Smillie.
\newblock Distribution of periodic points of polynomial diffeomorphisms of
  $\mathbf {C}\sp 2$.
\newblock {\em Invent. Math.}, 114(2):277--288, 1993.

\bibitem{BLS2}
E.~Bedford, M.~Lyubich, and J.~Smillie.
\newblock Polynomial diffeomorphisms of $\mathbf {C}\sp 2$. {I}{V}. {T}he
  measure of maximal entropy and laminar currents.
\newblock {\em Invent. Math.}, 112(1):77--125, 1993.

\bibitem{BS9}
E.~Bedford and J.~Smillie.
\newblock Real polynomial diffeomorphisms with maximal entropy: {T}angencies.
\newblock to appear, preprint available at [http://xxx.arxiv.org],
  arXiv:math.DS/0103038.

\bibitem{BS1}
E.~Bedford and J.~Smillie.
\newblock Polynomial diffeomorphisms of $\mathbf {C}\sp 2$: currents,
  equilibrium measure and hyperbolicity.
\newblock {\em Invent. Math.}, 103(1):69--99, 1991.

\bibitem{BS2}
E.~Bedford and J.~Smillie.
\newblock Polynomial diffeomorphisms of $\mathbf {C}\sp 2$. {I}{I}. {S}table
  manifolds and recurrence.
\newblock {\em J. Amer. Math. Soc.}, 4(4):657--679, 1991.

\bibitem{BS6}
E.~Bedford and J.~Smillie.
\newblock Polynomial diffeomorphisms of $\mathbf {C}\sp 2$. {V}{I}.
  {C}onnectivity of ${J}$.
\newblock {\em Ann. of Math. (2)}, 148(2):695--735, 1998.

\bibitem{BC}
M.~Benedicks and L.~Carleson.
\newblock The dynamics of the {H}\'enon map.
\newblock {\em Ann. of Math. (2)}, 133(1):73--169, 1991.

\bibitem{GenIA}
Interval Computations.
\newblock [http://www.cs.utep.edu/interval-comp/].

\bibitem{TCR}
T.~Cormen et~al.
\newblock {\em Introduction to Algorithms}.
\newblock The MIT Electrical Engineering and Computer Science Series. The MIT
  Press and McGraw-Hill Book Company, 1990.

\bibitem{Dell1}
M.~Dellnitz and O.~Junge.
\newblock Set oriented numerical methods for dynamical systems.
\newblock In {\em Handbook of dynamical systems, Vol. 2}, pages 221--264.
  North-Holland, Amsterdam, 2002.

\bibitem{DN}
R.~Devaney and Z.~Nitecki.
\newblock Shift automorphisms in the {H}\'enon mapping.
\newblock {\em Comm. Math. Phys.}, 67(2):137--146, 1979.

\bibitem{Eiden}
M.~Eidenschink.
\newblock {\em Exploring Global Dynamics: A Numerical Algorithm Based on the
  Conley Index Theory}.
\newblock PhD thesis, Georgia Institute of Technology, 1995.

\bibitem{FS1992c}
J.~E. Forn{\ae}ss and N.~Sibony.
\newblock Complex {H}\'enon mappings in $\mathbf {C}\sp 2$ and
  {F}atou-{B}ieberbach domains.
\newblock {\em Duke Math. J.}, 65(2):345--380, 1992.

\bibitem{FM}
S.~Friedland and J.~Milnor.
\newblock Dynamical properties of plane polynomial automorphisms.
\newblock {\em Ergodic Theory Dynamical Systems}, 9(1):67--99, 1989.

\bibitem{HK}
K.~Hoffman and R.~Kunze.
\newblock {\em Linear Algebra Second Edition}.
\newblock Prentice Hall, 1971.

\bibitem{SLHT}
J.~S.~L. Hruska.
\newblock {\em On the numerical construction of hyperbolic structures for
  complex dynamical systems}.
\newblock PhD thesis, Cornell University, 2002.
\newblock Download at [http://www.math.sunysb.edu/dynamics/theses/index.html].

\bibitem{SLHone}
S.~L. Hruska.
\newblock Rigorous numerical models for the dynamics of complex {H}\'{e}non
  mappings on their chain recurrent sets.
\newblock {\em Discrete and Continuous Dynamical Systems, to appear}.
\newblock available at [http://xxx.arxiv.org].

\bibitem{SLHtwo}
S.L. Hruska.
\newblock Constructing an expanding metric for dynamical systems in one complex
  variable.
\newblock {\em Nonlinearity}, 18:81--100, 2005.

\bibitem{HubKarl}
J.~Hubbard and K.~Papadantonakis.
\newblock Exploring the parameter space of complex {H}\'{e}non mappings.
\newblock {\em Journal of Experimental Mathematics}, to appear.

\bibitem{HPV}
J.~Hubbard, P.~Papadopol, and V.~Veselov.
\newblock A compactification of {H}\'{e}non mappings in $\mathbf {C}\sp 2$ as
  dynamical systems.
\newblock {\em Acta Math.}, 184(2):203--270, 2000.

\bibitem{HOV1}
J.~H. Hubbard and R.~W. Oberste-Vorth.
\newblock {H}\'{e}non mappings in the complex domain. {I}. {T}he global
  topology of dynamical space.
\newblock {\em Inst. Hautes \'Etudes Sci. Publ. Math.}, (79):5--46, 1994.

\bibitem{HOV2}
J.~H. Hubbard and R.~W. Oberste-Vorth.
\newblock {H}\'{e}non mappings in the complex domain. {I}{I}. projective and
  inductive limits of polynomials.
\newblock In B.~Branner and P.~Hjorth, editors, {\em Real and Complex Dynamical
  Systems}, volume 464 of {\em NATO Adv. Sci. Inst. Ser. C Math. Phys. Sci.},
  pages 89--132. Kluwer Acad. Publ., Dordrecht, 1995.

\bibitem{IS}
Y.~Ishii and J.~Smillie.
\newblock On the hyperbolicity of some complex {H}\'{e}non maps.
\newblock in preparation.

\bibitem{KMisch}
K.~Mischaikow.
\newblock Topological techniques for efficient rigorous computations in
  dynamics.
\newblock {\em Acta Numerica}, 2002.

\bibitem{MooreIA1}
R.~E. Moore.
\newblock {\em Interval Analysis}.
\newblock Prentice-Hall, Englewood Cliffs, New Jersey, 1966.

\bibitem{MooreIA2}
R.~E. Moore.
\newblock {\em Methods and Applications of Interval Analysis}.
\newblock SIAM Studies in Applied Mathematics, Philadelphia, 1979.

\bibitem{MNTU}
S.~Morosawa, Y.~Nishimura, M.~Taniguchi, and T.~Ueda.
\newblock {\em Holomorphic dynamics}.
\newblock Cambridge University Press, Cambridge, 2000.
\newblock Translated from the 1995 Japanese original and revised by the
  authors.

\bibitem{New}
S.~Newhouse.
\newblock Lectures on dynamical systems.
\newblock In {\em Dynamical Systems (C.I.M.E. Summer School, Bressanone,
  1978)}, volume~8 of {\em Progress in Mathematics}, pages 1--114.
  Birkh{\"{a}}user, Boston, Mass., 1980.

\bibitem{NP1}
S.~Newhouse and J.~Palis.
\newblock Bifurcations of {M}orse-{S}male dynamical systems.
\newblock In {\em Dynamical systems (Proc. Sympos., Univ. Bahia, Salvador,
  1971)}, pages 303--366. Academic Press, New York, 1973.

\bibitem{OVT}
R.~W. Oberste-Vorth.
\newblock {\em Complex horseshoes and the dynamics of mappings of two complex
  variables}.
\newblock PhD thesis, Cornell University, 1987.

\bibitem{OT}
R.~Oliva.
\newblock {\em On the combinatorics of external rays in the dynamics of the
  complex H\'{e}non map}.
\newblock PhD thesis, Cornell University, 1998.

\bibitem{Osi}
G.~Osipenko.
\newblock Construction of attractors and filtrations.
\newblock In {\em Conley index theory (Warsaw, 1997)}, volume~47 of {\em Banach
  Center Publ.}, pages 173--192. Polish Acad. Sci., Warsaw, 1999.

\bibitem{OsiCamp}
G.~Osipenko and S.~Campbell.
\newblock Applied symbolic dynamics: attractors and filtrations.
\newblock {\em Discrete Contin. Dynam. Systems}, 5(1):43--60, 1999.

\bibitem{PBIA}
PROFIL/BIAS Interval~Arithmetic Package.
\newblock \hfil\break
  [http://www.ti3.tu-harburg.de/Software/PROFILEnglisch.html].

\bibitem{Rob}
C.~Robinson.
\newblock {\em Dynamical systems}.
\newblock CRC Press, Boca Raton, FL, second edition, 1999.
\newblock Stability, symbolic dynamics, and chaos.

\bibitem{War}
W.~Tucker.
\newblock A rigorous {ODE} solver and {S}male's 14th problem.
\newblock {\em Found. Comput. Math.}, 2(1):53--117, 2002.

\end{thebibliography}

\end{document}